\documentclass[12pt]{amsart}
\usepackage{amscd, amsfonts, amssymb, amsthm}
\usepackage[all]{xypic}
\usepackage{graphics, lscape}
\usepackage{rotating}
  {\begin{list}
  {(\thesubsection.\theenumi)}
  {\usecounter{enumi}\setlength{\topsep}{8pt plus 2pt}
     \setlength{\parsep}{5.16663pt plus 1pt}
     \setlength{\itemsep}{5.16663pt plus 1pt}
     \setlength{\labelwidth}{2.5em}
     \setlength{\labelsep}{0.5em}
     \setlength{\itemindent}{0pt}
   }}
  {\end{list}}

\newcommand\CB{{\mathcal B}}
 
\renewcommand\CD{{\mathcal D}}
 
\newcommand\CF{{\mathcal F}}
\newcommand\CG{{\mathcal G}}
\newcommand\CH{{\mathcal H}}
\newcommand\CI{{\mathcal I}} 
\newcommand\CK{{\mathcal K}} 
\newcommand\CL{{\mathcal L}}
\newcommand\CM{{\mathcal M}} 
\newcommand\CN{{\mathcal N}}
\newcommand\CO{{\mathcal O}} 
\newcommand\CP{{\mathcal P}}
 
\newcommand\CQ{{\mathcal Q}}
 
\newcommand\CU{{\mathcal U}}

\newcommand\fb{{\mathfrak b}}
\newcommand\fg{{\mathfrak g}}
\newcommand\fk{{\mathfrak k}}
\newcommand\fp{{\mathfrak p}}
\newcommand\fq{{\mathfrak q}}
\newcommand\fu{{\mathfrak u}} 
\newcommand\fl{{\mathfrak l}}

\newcommand\ft{{\mathfrak t}}

\newcommand\FD{\mathfrak D}
\newcommand\FN{\mathfrak N}
\newcommand\FO{{\mathfrak O}}
\newcommand\FS{\mathfrak S}
\newcommand\FV{\mathfrak V}

\newcommand\BBC{{\mathbb C}}
\newcommand\BBD{{\mathbb D}}
\newcommand\BBQ{{\mathbb Q}}
\newcommand\BBR{{\mathbb R}}
\newcommand\BBZ{{\mathbb Z}}

\newcommand\Coinv{\operatorname{Coinv}}
\newcommand\Ext{{\operatorname{Ext}}}
\newcommand\End{{\operatorname{End}}}
\newcommand\Hom{{\operatorname{Hom}}}
\newcommand\GL{\operatorname{GL}}
\newcommand\gr{{\operatorname{gr}}}
\newcommand\IC{{\operatorname{IC}}}
\newcommand\Ind{{\operatorname{Ind}}}
\newcommand\Lie{{\operatorname{Lie}}}
\newcommand\Mod{{\operatorname{Mod}}}

\newcommand\SL{\operatorname{SL}}

\newcommand\inverse{^{-1}}
\renewcommand\th{{^{\text{th}}}}

\newcommand\alphacheck{\check \alpha}

\newcommand\coh{\text{coh}}
\newcommand\FObar{\overline{\mathfrak O}}
\newcommand\FNt{{\widetilde{\FN}}}
\newcommand\fgt{{\widetilde{\fg}}}

\newcommand\fgrst{{\widetilde {\fg}_{\rs}}}
\newcommand\fingen{\text{fg}}
\newcommand\FVbar{\overline{\FV}}
\newcommand\Gbar{\overline{G}}
\newcommand\id{{id}}
\newcommand\op{\text{op}}

\newcommand\reg{{\operatorname{reg}}}
\newcommand\rs{{\operatorname{rs}}}
\newcommand\sbar{\overline{\langle s \rangle}} 
\newcommand\Tbar{\overline{T}}
\newcommand\wtilde{{\widetilde{w}}} 
\newcommand\XPQ{X^{\CP, \CQ}}
\newcommand\XPQcd{X^{\CP, \CQ}_{c,d}} 
\newcommand\XPQoo{X^{\CP,\CQ}_{0,0}} 
\newcommand\XPQrr{X^{\CP, \CQ}_{\reg,\reg}}
\newcommand\Ywbar{{\overline{Y_w}}} 
\newcommand\Zhat{{\widehat Z}}
\newcommand\ZPQ{Z^{\CP, \CQ}} 
\newcommand\ZPQcd{Z^{\CP, \CQ}_{c,d}}
\newcommand\ZPQw{Z^{\CP, \CQ}_w} 
\newcommand\Zwbar{{\overline{Z_w}}}

\numberwithin{equation}{section}

\theoremstyle{plain}
\newtheorem{lemma}[equation]{Lemma}
\newtheorem{theorem}[equation]{Theorem}
\newtheorem{conjecture}[equation]{Conjecture}
\newtheorem{corollary}[equation]{Corollary}
\newtheorem{proposition}[equation]{Proposition}

\thanks{The authors would like to thank their charming wives for their
  unwavering support during the preparation of this paper}

\subjclass[2000]{Primary 22E46, 19L47, 20G05; Secondary 14F99, 20G99}

\begin{document}

\title[Steinberg variety and representations] {The Steinberg variety\\
  and representations of reductive groups}

\dedicatory{Dedicated to Gus Lehrer on the occasion of his 60th birthday}

\author[J.M. Douglass]{J. Matthew Douglass} \address{Department of
  Mathematics\\ University of North Texas\\ Denton TX, USA 76203}
\email{douglass@unt.edu} \urladdr{http://hilbert.math.unt.edu}

\author[G. R\"ohrle]{Gerhard R\"ohrle} \address
{Fakult\"at f\"ur Mathematik, Ruhr-Universit\"at Bochum, D-44780
  Bochum, Germany} \email{gerhard.roehrle@rub.de}
\urladdr{http://www.ruhr-uni-bochum.de/ffm/Lehrstuehle/Lehrstuhl-VI/rubroehrle.html}
 
\date{\today}

\maketitle
\allowdisplaybreaks

\begin{abstract} 
  We give an overview of some of the main results in geometric
  representation theory that have been proved by means of the
  Steinberg variety.  Steinberg's insight was to use such a variety of
  triples in order to prove a conjectured formula by Grothendieck.
  The Steinberg variety was later used to give an alternative approach
  to Springer's representations and played a central role in the proof
  of the Deligne-Langlands conjecture for Hecke algebras by Kazhdan
  and Lusztig.
\end{abstract}


\section{Introduction}

Suppose $G$ is a connected, reductive algebraic group defined over an
algebraically closed field $k$, $\CB$ is the variety of Borel
subgroups of $G$, and $u$ is a unipotent element in $G$. Let $\CB_u$
denote the closed subvariety of $\CB$ consisting of those Borel
subgroups that contain $u$, let $r$ denote the rank of $G$, and let
$C$ denote the conjugacy class of $u$.

In 1976, motivated by the problem of proving the equality conjectured
by Grothendieck
\begin{equation}\tag{$*$}
  \label{*}
  \dim Z_G(u)= r+ 2 \dim \CB_u,
\end{equation}
in order to get the multiplicity $2$ in (\ref{*}) in the picture,
Steinberg \cite{steinberg:desingularization} introduced a variety of
triples
\[
S=\{\, (v, B, B')\in C\times \CB\times \CB\mid v\in B\cap B'\,\}.
\]
By analyzing the geometry of the variety $S$, he was able to prove
(\ref{*}) in most cases. In addition, by exploiting the fact that the
$G$-orbits on $\CB\times \CB$ are canonically indexed by elements of
the Weyl group of $G$, he showed that $S$ could be used to establish
relationships between Weyl group elements and unipotent elements in
$G$.

Now let $\fg$ denote the Lie algebra of $G$, and let $\FN$ denote the
variety of nilpotent elements in $\fg$. The \emph{Steinberg variety of
  $G$} is
\[
Z=\{\, (x, B, B')\in \FN\times \CB\times \CB\mid x\in \Lie(B)\cap
\Lie(B')\,\}.
\]
If the characteristic of $k$ is zero or good for $G$, then there is a
$G$-equivariant isomorphism between $\FN$ and $\CU$, the variety of
unipotent elements in $G$, and so $Z\cong \{\, (u, B, B')\in \CU
\times \CB\times \CB\mid u\in B\cap B'\,\}$.

In the thirty years since Steinberg first exploited the variety $S$,
the Steinberg variety has played a key role in advancing our
understanding of objects that at first seem to be quite unrelated:
\begin{itemize}
\item Representations of the Weyl group $W$ of $G$
\item The geometry of nilpotent orbits in $\fg$ and their covers
\item Differential operators on $\CB$
\item Primitive ideals in the universal enveloping algebra of $\fg$
\item Representations of $p$-adic groups and the local Langlands
  program
\end{itemize}

In this paper we hope to give readers who are familiar with some
aspects of the representation theory of semisimple algebraic groups,
or Lie groups, but who are not specialists in this particular flavor
of geometric representation theory, an overview of the main results
that have been proved using the Steinberg variety. In the process we
hope to make these results more accessible to non-experts and at the
same time emphasize the unifying role played by the Steinberg
variety. We assume that the reader is quite familiar with the basics
of the study of algebraic groups, especially reductive algebraic
groups and their Lie algebras, as contained in the books by Springer
\cite{springer:algebraic} and Carter \cite{carter:finite} for example.
 
We will more or less follow the historical development, beginning with
concrete, geometric constructions and then progressing to increasingly
more advanced and abstract notions. 

In \S2 we analyze the geometry of $Z$, including applications to
orbital varieties, characteristic varieties and primitive ideals, and
generalizations.

In \S3 we study the Borel-Moore homology of $Z$ and the relation with
representations of Weyl groups. Soon after Steinberg introduced his
variety $S$, Kazhdan and Lusztig \cite{kazhdanlusztig:topological}
defined an action of $W\times W$ on the top Borel-Moore homology group
of $Z$. Following a suggestion of Springer, they showed that the
representation of $W\times W$ on the top homology group, $H_{4n}(Z)$,
is the two-sided regular representation of $W$. Somewhat later,
Ginzburg \cite{ginzburg:G-modules} and independently Kashiwara and
Tanisaki \cite{kashiwaratanisaki:characteristic}, defined a
multiplication on the total Borel-Moore homology of $Z$. With this
multiplication, $H_{4n}(Z)$ is a subalgebra isomorphic to the group
algebra of $W$.

The authors \cite{douglassroehrle:homology}
\cite{douglassroehrle:coinvariant} have used Ginzburg's construction
to describe the top Borel-Moore homology groups of the generalized
Steinberg varieties $\XPQoo$ and $\XPQrr$ (see \S\ref{s2.4}) in terms
of $W$, as well as to give an explicit, elementary, computation of the
total Borel-Moore homology of $Z$ as a graded algebra: it is
isomorphic to the smash product of the coinvariant algebra of $W$ and
the group algebra of $W$.

Orbital varieties arise naturally in the geometry of the Steinberg
variety. Using the convolution product formalism, Hinich and Joseph
\cite{hinichjoseph:orbital} have recently proved an old conjecture of
Joseph about inclusions of closures of orbital varieties.

In \S4 we study the equivariant $K$-theory of $Z$ and what is
undoubtedly the most important result to date involving the Steinberg
variety: the Kazhdan-Lusztig isomorphism
\cite{kazhdanlusztig:langlands} between $K^{G\times \BBC^*}(Z)$ and
the extended, affine Hecke algebra $\CH$. Using this isomorphism,
Kazhdan and Lusztig were able to classify the irreducible
representations of $\CH$ and hence to classify the representations
containing a vector fixed by an Iwahori subgroup of the $p$-adic group
with the same type as the Langlands dual ${}^LG$ of $G$. In this way,
the Steinberg variety plays a key role in the local Langlands program
and also leads to a better understanding of the extended affine Hecke
algebra.

Very recent work involving the Steinberg variety centers around
attempts to categorify the isomorphism between the specialization of
$K^{G\times \BBC^*}(Z)$ at $p$ and the Hecke algebra of Iwahori
bi-invariant functions on ${}^LG(\BBQ_p)$. Because of time and space
constraints, we leave a discussion of this research to a future
article.

\section{Geometry}

For the rest of this paper, in order to simplify the exposition, we
assume that $G$ is connected, the derived group of $G$ is simply
connected, and that $k=\BBC$. Most of the results below hold, with
obvious modifications, for an arbitrary reductive algebraic group when
the characteristic of $k$ is zero or very good for $G$ (for the
definition of ``very good characteristic'' see \cite[\S1.14]
{carter:finite}).

Fix a Borel subgroup $B$ in $G$ and a maximal torus $T$ in $B$. Define
$U$ to be the unipotent radical of $B$ and define $W=N_G(T)/T$ to be
the Weyl group of $(G,T)$. Set $n= \dim \CB$ and $r=\dim T$.

We will use the convention that a lowercase fraktur letter denotes the
Lie algebra of the algebraic group denoted by the corresponding
uppercase roman letter.

For $x$ in $\FN$, define $\CB_x= \{\, gBg\inverse\mid g\inverse x\in
\fb\,\}$, the \emph{Springer fibre} at $x$.

\subsection{Irreducible components of $Z$, Weyl group elements, and
  nilpotent orbits} \label{s2.1}

We begin analyzing the geometry of $Z$ using ideas that go back to
Steinberg \cite{steinberg:desingularization} and Spaltenstein
\cite{spaltenstein:classes}.

The group $G$ acts on $\CB$ by conjugation and on $\FN$ by the adjoint
action. This latter action is denoted by $(g,x)\mapsto g\cdot
x=gx$. Thus, $G$ acts ``diagonally'' on $Z$.

Let $\pi\colon Z\to \CB\times \CB$ be the projection on the second and
third factors. By the Bruhat Lemma, the elements of $W$ parametrize
the $G$-orbits on $\CB\times \CB$. An element $w$ in $W$ corresponds
to the $G$-orbit containing $(B, wBw\inverse)$ in $\CB \times \CB$.
Define
\[
Z_w= \pi\inverse \left(G(B, wBw\inverse) \right),\ U_w= U\cap
wUw\inverse,\ \text{and}\ B_w= B\cap wBw\inverse.
\]
The varieties $Z_w$ play a key role in the rest of this paper.

For $w$ in $W$, the restriction of $\pi$ to $Z_w$ is a $G$-equivariant
morphism from $Z_w$ onto a transitive $G$-space. The fibre over the
point $(B, wBw\inverse)$ is isomorphic to $\fu_w$ and so it follows
from \cite[II 3.7]{slodowy:simple} that $Z_w$ is isomorphic to the
associated fibre bundle $G\times^{B_w} \fu_w$.  Thus, $Z_w$ is
irreducible and $\dim Z_w= \dim G-\dim B_w+\dim \fu_w= 2n$.
Furthermore, each $Z_w$ is locally closed in $Z$ and so it follows
that $\{\,\Zwbar\mid w\in W\,\}$ is the set of irreducible components
of $Z$.

Now let $\mu_z\colon Z\to \FN$ denote the projection on the first
component.  For a $G$-orbit, $\FO$, in $\FN$, set $Z_\FO=
\mu_z\inverse (\FO)$ and fix $x$ in $\FO$. Then the restriction of
$\mu_z$ to $Z_{\FO}$ is a $G$-equivariant morphism from $Z_\FO$ onto a
transitive $G$-space. The fibre over $x$ is isomorphic to $\CB_x\times
\CB_x$ and so it follows from \cite[II 3.7]{slodowy:simple} that
$Z_\FO\cong G\times^{Z_G(x)} (\CB_x\times \CB_x)$. Spaltenstein
\cite[\S II.1]{spaltenstein:classes} has shown that the variety
$\CB_x$ is equidimensional and Steinberg and Spaltenstein have shown
that $\dim Z_G(x)= r+2\dim \CB_x$. This implies the following results
due to Steinberg \cite[Proposition 3.1]{steinberg:desingularization}:
\begin{itemize}
\item[(1)] $\dim Z_\FO= \dim G-\dim Z_G(x)+ 2\dim \CB_x = \dim
  G-r=2n$.
\item[(2)] Every irreducible component of $Z_\FO$ has the form
  \[
  G (\{x\} \times C_1\times C_2)= G(\{x\} \times (Z_G(x)(C_1\times
  C_2)))
  \]
  where $C_1$ and $C_2$ are irreducible components of $\CB_x$.
\item[(3)] A pair, $(C_1',C_2')$, of irreducible components of $\CB_x$
  determines the same irreducible component of $Z_\FO$ as $(C_1, C_2)$
  if and only if there is a $z$ in $Z_G(x)$ with $(C_1',C_2') =(zC_1
  z\inverse ,zC_2z\inverse)$.
\end{itemize}

From (2) we see that $Z_\FO$ is equidimensional with $\dim Z_\FO=
2n=\dim Z$ and from (3) we see that there is a bijection between
irreducible components of $Z_\FO$ and $Z_G(x)$-orbits on the set of
irreducible components of $\CB_x \times \CB_x$.

The closures of the irreducible components of $Z_{\FO}$ are closed,
irreducible, $2n$-dimensional subvarieties of $Z$ and so each
irreducible component of $Z_{\FO}$ is of the form $Z_{\FO} \cap
\Zwbar$ for some unique $w$ in $W$. Define $W_{\FO}$ to be the subset
of $W$ that parametrizes the irreducible components of $Z_{\FO}$.
Then $w$ is in $W_{\FO}$ if and only if $Z_\FO\cap \Zwbar$ is an
irreducible component of $Z_{\FO}$.

Clearly, $W$ is the disjoint union of the $W_\FO$'s as $\FO$ varies
over the nilpotent orbits in $\FN$. The subsets $W_\FO$ are called
\emph{two-sided Steinberg cells}. Two-sided Steinberg cells have
several properties in common with two-sided Kazhdan-Lusztig cells in
$W$. Some of the properties of two-sided Steinberg cells will be
described in the next subsection.  Kazhdan-Lusztig cells were
introduced in \cite[\S 1]{kazhdanlusztig:coxeter}. We will briefly
review this theory in \S\ref{4.4}.

In general there are more two-sided Steinberg cells than two-sided
Kazhdan-Lusztig cells. This may be seen as follows. Clearly, two-sided
Steinberg cells are in bijection with the set of $G$-orbits in $\FN$.

Two-sided Kazhdan-Lusztig cells may be related to nilpotent orbits
through the Springer correspondence using Lusztig's analysis of
Kazhdan-Lusztig cells in Weyl groups. We will review the Springer
correspondence in \S\ref{s3.4} below, where we will see that there is
an injection from the set of nilpotent orbits to the set of
irreducible representations of $W$ given by associating with $\FO$ the
representation of $W$ on $H_{2d_x}(\CB_x)^{C(x)}$, where $x$ is in
$\FO$ and $C(x)$ is the component group of $x$. Two-sided
Kazhdan-Lusztig cells determine a filtration of the group algebra
$\BBQ[W]$ by two-sided ideals (see \S\ref{4.4}) and in the associated
graded $W\times W$-module, each summand contains a distinguished
representation that is called \emph{special} (see \cite{lusztig:class}
and \cite[Chapter 5]{lusztig:characters}). The case-by-case
computation of the Springer correspondence shows that every special
representation of $W$ is equivalent to the representation of $W$ on
$H_{2d_x}(\CB_x)^{C(x)}$ for some $x$. The resulting nilpotent orbits
are called \emph{special} nilpotent orbits.

If $G$ has type $A_l$, then every irreducible representation of $W$
and every nilpotent orbit is special but otherwise there are
non-special irreducible representation of $W$ and nilpotent orbits.
Although in general there are fewer two-sided Kazhdan-Lusztig cells in
$W$ than two-sided Steinberg cells, Lusztig \cite[\S
4]{lusztig:cellsIV} has constructed a bijection between the set of
two-sided Kazhdan-Lusztig cells in the extended, affine, Weyl group,
$W_e$, and the set of $G$-orbits in $\FN$. Thus, there is a bijection
between two-sided Steinberg cells in $W$ and two-sided Kazhdan-Lusztig
cells in $W_e$. We will describe this bijection in \S\ref{4.4} in
connection with the computation of the equivariant $K$-theory of the
Steinberg variety.

Suppose $\FO$ is a nilpotent orbit and $x$ is in $\FO$. We can
explicitly describe the bijection in (c) above between $W_\FO$ and the
$Z_G(x)$-orbits on the set of pairs of irreducible components of
$\CB_x$ as follows. If $w$ is in $W_{\FO}$ and $(C_1,C_2)$ is a pair
of irreducible components of $\CB_x$, then $w$ corresponds to the
$Z_G(x)$-orbit of $(C_1, C_2)$ if and only if $G(B, wBw\inverse) \cap
(C_1\times C_2)$ is dense in $C_1\times C_2$.

Using the isomorphism $Z_w\cong G\times^{B_w} \fu_w$ we see that
$Z_\FO \cap Z_w\cong G\times^{B_w} (\FO \cap \fu_w)$. Therefore, $w$
is in $W_\FO$ if and only if $\FO\cap \fu_w$ is dense in $\fu_w$. This
shows in particular that $W_\FO$ is closed under taking inverses.

We conclude this subsection with some examples of two-sided Steinberg
cells.

When $x=0$ we have $Z_{\{0\}}= \overline{Z_{w_0}}= \{0\}\times
\CB\times \CB$ where $w_0$ is the longest element in $W$. Therefore,
$W_{\{0\}}= \{w_0\}$. 

At the other extreme, let $\FN_\reg$ denote the regular nilpotent
orbit. Then it follows from the fact that every regular nilpotent
element is contained in a unique Borel subalgebra that $W_{\FN_\reg}$
contains just the identity element in $W$.

For $G$ of type $A_l$, it follows from a result of Spaltenstein
\cite{spaltenstein:fixed} that two elements of $W$ lie in the same
two-sided Steinberg cell if and only if they yield the same Young
diagram under the Robinson-Schensted correspondence. A more refined
result due to Steinberg will be discussed at the end of the next
subsection.

\subsection{Orbital varieties}\label{s2.2}

Suppose that $\FO$ is a nilpotent orbit. An \emph{orbital variety for
  $\FO$} is an irreducible component of $\FO\cap \fu$. An
\emph{orbital variety} is a subvariety of $\FN$ that is orbital for
some nilpotent orbit. The reader should be aware that sometimes an
orbital variety is defined as the closure of an irreducible component
of $\FO\cap \fu$.

We will see in this subsection that orbital varieties can be used to
decompose two-sided Steinberg cells into left and right Steinberg
cells and to refine the relationship between nilpotent orbits and
elements of $W$. When $G$ is of type $A_l$ and $W$ is the symmetric
group $S_{l+1}$, the decomposition of a two-sided Steinberg cell into
left and right Steinberg cells can be viewed as a geometric
realization of the Robinson-Schensted correspondence.

We will see in the next subsection that orbital varieties arise in the
theory of associated varieties of finitely generated $\fg$-modules.

Fix a nilpotent orbit $\FO$ and an element $x$ in $\FO\cap \fu$.
Define $p\colon G\to \FO$ by $p(g)= g\inverse x$ and $q\colon G\to
\CB$ by $q(g)= gBg\inverse$.  Then $p\inverse( \FO\cap \fu)=
q\inverse( \CB_x)$. Spaltenstein \cite[\S II.2]{spaltenstein:classes}
has shown that
\begin{itemize}
\item[(1)] if $C$ is an irreducible component of $\CB_x$, then
  $pq\inverse(C)$ is an orbital variety for $\FO$,
\item[(2)] every orbital variety for $\FO$ has the form
  $pq\inverse(C)$ for some irreducible component $C$ of $\CB_x$, and
\item[(3)] $pq\inverse(C) =pq \inverse(C')$ for components $C$ and
  $C'$ of $\CB_x$ if and only if $C$ and $C'$ are in the same
  $Z_G(x)$-orbit.
\end{itemize} 
It follows immediately that $\FO\cap \fu$ is equidimensional and all
orbital varieties for $\FO$ have the same dimension: $n-\dim
\CB_x=\frac 12 \dim \FO$.

We decompose two-sided Steinberg cells into left and right
Steinberg cells following a construction of Joseph
\cite[\S9]{joseph:variety}.

Suppose $\FV_1$ and $\FV_2$ are orbital varieties for $\FO$. Choose
irreducible components $C_1$ and $C_2$ of $\CB_x$ so that
$pq\inverse(C_1)= \FV_1$ and $pq\inverse(C_2)= \FV_2$. We have seen
that there is a $w$ in $W_\FO$ so that $Z_\FO \cap \Zwbar= G \left(
  \{x\}\times Z_G(x)(C_1 \times C_2) \right)$. Clearly,
$\overline{\mu_z\inverse(x) \cap Z_w} \subseteq \mu_z\inverse(x) \cap
\Zwbar$.  Since both sides are closed, both sides are $Z_G(x)$-stable,
and the right hand side is the $Z_G(x)$-saturation of $\{x\} \times
C_1 \times C_2$, it follows that $\overline{ \mu_z\inverse(x) \cap
  Z_w} = \mu_z\inverse(x) \cap \Zwbar$.

Let $p_2$ denote the projection of $Z_\FO$ to $\CB$ by $p_2(x, B',
B'')= B'$. Then $pq\inverse p_2\left( \mu_z\inverse(x) \cap Z_w
\right)= B(\FO\cap \fu_w)$. Also,
\[
pq\inverse p_2\left(\mu_z\inverse(x) \cap \Zwbar \right) =pq\inverse
p_2\left( \{x\}\times Z_G(x)(C_1 \times C_2) \right) =
pq\inverse\left( Z_G(x) C_1 \right) =\FV_1.
\]
Since $\FO\cap \fu_w$ is dense in $\fu_w$ we have $\overline{ B\fu_w}
\cap \FO= \overline{B(\FO\cap \fu_w)} \subseteq \FV_1$. However, since
$\mu_z\inverse(x) \cap Z_w$ is a dense, $Z_G(x)$-stable subset of
$\mu_z\inverse(x)\cap \Zwbar$, it follows that
\begin{align*}
  \dim B(\FO\cap \fu_w)&= \dim pq\inverse p_2\left( \mu_z\inverse(x)
    \cap Z_w \right)\\
  &= \dim p_2\left(\mu_z\inverse(x) \cap \Zwbar \right) +\dim B
  -\dim Z_G(x)\\
  &= \dim \CB_x+\dim B-r-2\dim \CB_x\\
  &= n-\dim \CB_x
\end{align*}
and so $\overline{ B\fu_w} \cap \FO= \FV_1$.

A similar argument shows that $\overline{ B\fu_{w\inverse}} \cap \FO=
\FV_2$. This proves the following theorem.

\begin{theorem}\label{pairs}
  If $\FO$ is a nilpotent orbit and $\FV_1$ and $\FV_2$ are orbital
  varieties for $\FO$, then there is a $w$ in $W_\FO$ so that $\FV_1=
  \overline{B \fu_w} \cap \FO$ and $\FV_2= \overline{B
    \fu_{w\inverse}} \cap \FO$.
\end{theorem}

Conversely, if $w$ is in $W_\FO$, then $\fu_w$ is irreducible and the
arguments above show that $\fu_w\cap \FO$ is dense in $\fu_w$ and then
that $\overline{B \fu_w} \cap \FO$ is an orbital variety. This proves
the next proposition.

\begin{proposition}
  Orbital varieties are the subsets of $\fu$ of the form $\overline{B
    \fu_w} \cap \FO$, where $\fu_w\cap \FO$ is dense in $\fu_w$.
\end{proposition}

For $w$ in $W$, define $\FV_l(w)= \overline{ B\fu_{w\inverse}} \cap
\FO$ when $w$ is in $W_\FO$. For $w_1$ and $w_2$ in $W$, define
$w_1\sim_l w_2$ if $\FV_l(w_1)= \FV_l(w_2)$. Then $\sim_l$ is an
equivalence relation and the equivalence classes are called \emph{left
  Steinberg cells.}  Similarly, define $\FV_r(w)= \overline{ B\fu_w}
\cap \FO$ when $w$ is in $W_\FO$ and $w_1\sim_r w_2$ if $\FV_r(w_1)=
\FV_r(w_2)$. The equivalence classes for $\sim_r$ are called
\emph{right Steinberg cells.}

Clearly, each two-sided Steinberg cell is a disjoint union of left
Steinberg cells and is also the disjoint union of right Steinberg
cells.  Precisely, if $w$ is in $W_\FO$, then
\[
W_\FO= \coprod_{y\in \FV_r(w)} \FV_l(y)= \coprod_{y\in \FV_l(w)}
\FV_r(y).
\]

It follows from Theorem \ref{pairs} that the rule $w\mapsto (\FV_r(w),
\FV_l(w))$ defines a surjection from $W$ to the set of pairs of
orbital varieties for the same nilpotent orbit. We will see in
\S\ref{s3.4} that the number of orbital varieties for a nilpotent
orbit $\FO$ is the dimension of the Springer representation of $W$
corresponding to the trivial representation of the component group of
any element in $\FO$. Denote this representation of $W$ by $\rho_\FO$.
Then the number of pairs $(\FV_1, \FV_2)$, where $\FV_1$ and $\FV_2$
are orbital varieties for the same nilpotent orbit, is $\sum_{\FO}
(\dim \rho_\FO)^2$. In general this sum is strictly smaller than
$|W|$. Equivalently, in general, there are more irreducible
representations of $W$ than $G$-orbits in $\FN$.

However, if $G$ has type $A$, for example if $G=\SL_n(\BBC)$ or
$\GL_n(\BBC)$, then every irreducible representation of $W$ is of the
form $\rho_\FO$ for a unique nilpotent orbit $\FO$. In this case
$w\mapsto (\FV_r(w), \FV_l(w))$ defines a bijection from $W$ to the
set of pairs of orbital varieties for the same nilpotent orbit.
Steinberg has shown that this bijection is essentially given by the
Robinson-Schensted correspondence.

In more detail, using the notation in the proof of Theorem
\ref{pairs}, suppose that $\FO$ is a nilpotent orbit, $\FV_1$ and
$\FV_2$ are orbital varieties for $\FO$, and $C_1$ and $C_2$ are the
corresponding irreducible components in $\CB_x$. In
\cite{steinberg:robinson} Steinberg defines a function from $\CB$ to
the set of standard Young tableaux and shows that $G(B, wBw\inverse)
\cap (C_1\times C_2)$ is dense in $C_1\times C_2$ if and only if the
pair of standard Young tableaux associated to a generic pair $(B',
B'')$ in $C_1\times C_2$ is the same as the pair of standard Young
tableaux associated to $w$ by the Robinson-Schensted correspondence.
For more details, see also \cite{douglass:involution}.

An open problem, even in type $A$, is determining the orbit closures
of orbital varieties. Some rudimentary information may be obtained by
considering the top Borel-Moore homology group of $Z$ (see \S3 below
and \cite[\S4, \S5]{hinichjoseph:orbital}).

\subsection{Associated varieties and characteristic varieties}
\label{s2.3}

The Steinberg variety and orbital varieties also arise naturally in
the Beilinson-Bernstein theory of algebraic $(\CD, K)$-modules
\cite{beilinsonbernstein:localisation}. This was first observed by
Borho and Brylinski \cite{borhobrylinski:differentialIII} and Ginzburg
\cite{ginzburg:G-modules}. In this subsection we begin with a review
of the Beilinson-Bernstein Localization Theorem and its connection
with the computation of characteristic varieties and associated
varieties. Then we describe an equivariant version of this theory. It
is in the equivariant theory that the Steinberg variety naturally
occurs.

For a variety $X$ (over $\BBC$), let $\CO_X$ denote the structure
sheaf of $X$, $\BBC[X]=\Gamma(X, \CO_X)$ the algebra of global,
regular functions on $X$, and $\CD_X$ the sheaf of algebraic
differential operators on $X$. On an open subvariety, $V$, of $X$,
$\Gamma(V, \CD_X)$ is the subalgebra of $\Hom_{\BBC}(\BBC[V],
\BBC[V])$ generated by multiplication by elements of $\BBC[V]$ and
$\BBC$-linear derivations of $\BBC[V]$. Define $D_X= \Gamma(X,
\CD_X)$, the algebra of global, algebraic, differential operators on
$X$.

A \emph{quasi-coherent $\CD_X$-module} is a left $\CD_X$-module that
is quasi-coherent when considered as an $\CO_X$-module. Generalizing a
familiar result for affine varieties, Beilinson-Bernstein
\cite[\S2]{beilinsonbernstein:localisation} have proved that for
$X=\CB$, the global section functor, $\Gamma(\CB, \,\cdot\,)$, defines
an equivalence of categories between the category of quasi-coherent
$\CD_\CB$-modules and the category of $D_\CB$-modules.

In turn, the algebra $D_\CB$ is isomorphic to $U(\fg)/I_0$, where
$U(\fg)$ is the universal enveloping algebra of $\fg$ and $I_0$
denotes the two-sided ideal in $U(\fg)$ generated by the kernel of the
trivial character of the center of $U(\fg)$ (see \cite[\S
3]{borhobrylinski:differentialI}). Thus, the category of
$D_\CB$-modules is equivalent to the category of $U(\fg)/I_0$-modules,
that is, the category of $U(\fg)$-modules with trivial central
character.

Composing these two equivalences we see that the category of
quasi-coherent $\CD_\CB$-modules is equivalent to the category of
$U(\fg)$-modules with trivial central character. In this equivalence,
coherent $\CD_\CB$-modules (that is, $\CD_\CB$-modules that are
coherent when considered as $\CO_\CB$-modules) correspond to finitely
generated $U(\fg)$-modules with trivial central character.

The equivalence of categories between coherent $\CD_\CB$-modules and
finitely generated $U(\fg)$-modules with trivial central character has
a geometric shadow that can be described using the ``moment map'' of
the $G$-action on the cotangent bundle of $\CB$.

Let $B'$ be a Borel subgroup of $G$. Then using the Killing form on
$\fg$, the cotangent space to $\CB$ at $B'$ may be identified with
$\fb'\cap \FN$, the nilradical of $\fb'$.  Define
\[
\FNt=\{\,(x, B')\in \FN\times \CB\mid x\in \fb'\,\}
\]
and let $\mu\colon \FNt\to \FN$ be the projection on the first factor.
Then $\FNt\cong T^*\CB$, the cotangent bundle of $\CB$. It is easy to
see that $Z\cong \FNt \times_\FN \FNt \cong T^*\CB \times_\FN T^*
\CB$.

Using the orders of differential operators, we obtain a filtration of
$\CD_X$. With respect to this filtration, the associated graded sheaf
$\gr\, \CD_\CB$ is isomorphic to the direct image $p_* \CO_{T^*\CB}$,
where $p\colon T^*\CB \to \CB$ is the projection. 

Let $\CM$ be a coherent $\CD_\CB$-module. Then $\CM$ has a ``good''
filtration such that $\gr\, \CM$ is a coherent $\gr\, \CD_\CB$-module.
Since $\gr\, \CD_\CB \cong p_* \CO_{T^*\CB}$, we see that $\gr\, \CM$
has the structure of a coherent $\CO_{T^*\CB}$-module. The
\emph{characteristic variety of $\CM$} is the support in $T^*\CB$ of
the $\CO_{T^*\CB}$-module $\gr\, \CM$. Using the isomorphism $T^*\CB
\cong \FNt$, we identify the characteristic variety of $\CM$ with a
closed subvariety of $\FNt$ and denote this latter variety by
$V_\FNt(\CM)$. It is known that $V_{\FNt}(\CM)$ is independent of the
choice of good filtration.

Now consider the enveloping algebra $U(\fg)$ with the standard
filtration. By the PBW Theorem, $\gr\,U(\fg) \cong \operatorname{
  Sym}(\fg)$, the symmetric algebra of $\fg$. Using the Killing form,
we identify $\fg$ with its linear dual, $\fg^*$, and $\gr\, U(\fg)$
with $\BBC[\fg]$. Let $M$ be a finitely generated $U(\fg)$-module.
Then $M$ has a ``good'' filtration such that the associated graded
module, $\gr\, M$, a module for $\gr\, U(\fg)\cong \BBC[\fg]$, is
finitely generated. The \emph{associated variety of $M$}, denoted by
$V_{\fg}(M)$, is the support of the $\BBC[\fg]$-module $\gr\, M$ -- a
closed subvariety of $\fg$. It is known that $V_{\fg}(M)$ is
independent of the choice of good filtration.

Borho and Brylinski \cite[\S1.9]{borhobrylinski:differentialIII} have
proved the following theorem.

\begin{theorem}\label{thmsigmaa}
  Suppose that $\CM$ is a coherent $\CD_\CB$-module and let $M$ denote
  the space of global sections of $\CM$. Then $V_\fg(M)\subseteq \FN$
  and $\mu( V_{\FNt}(\CM))= V_\fg(M)$.
\end{theorem}

There are equivariant versions of the above constructions which
incorporate a subgroup of $G$ that acts on $\CB$ with finitely many
orbits. It is in this equivariant context that the Steinberg variety
and orbital varieties make their appearance.

Suppose that $K$ is a closed, connected, algebraic subgroup of $G$
that acts on $\CB$ with finitely many orbits. The two special cases we
are interested in are the ``highest weight'' case, when $K=B$ is a
Borel subgroup of $G$, and the ``Harish-Chandra'' case, when $K=G_d$
is the diagonal subgroup of $G\times G$.

In the general setting, we suppose that $W$ is a finite set that
indexes the $K$-orbits on $\CB$ by $w \leftrightarrow X_w$.  Of
course, in the examples we are interested in, we know that the Weyl
group $W$ indexes the set of orbits of $K$ on $\CB$.

For $w$ in $W$, let $T^*_{w} \CB$ denote the conormal bundle to the
$K$-orbit $X_w$ in $T^*\CB$.  Then letting $\fk^\perp$ denote the
subspace of $\fg$ orthogonal to $\fk$ with respect to the Killing form
and using our identification of $T^*\CB$ with pairs, we may identify
\[
T^*_{w} \CB= \{\, (x, B') \in \FN\times \CB\mid B'\in X_w,\ x\in
\fb' \cap \fk^\perp\,\}.
\]
Define $Y_{\fk^\perp}= \mu \inverse (\fk^\perp \cap \FN)$. Then
$Y_{\fk^\perp}$ is closed, $Y_{\fk^\perp} =\coprod_{w\in W} T^*_{w}
\CB =\cup_{w\in W} \overline{T_w^* \CB}$, and $\mu$ restricts to a
surjection $Y_{\fk^\perp} \xrightarrow{\mu} \fk^\perp$ (see \cite[\S
2.4]{borhobrylinski:differentialIII}).  Summarizing, we have a
commutative diagram
\begin{equation}
  \label{diagY}
  \xymatrix{Y_{\fk^\perp} \ar[r] \ar[d]_{\mu} &\FNt \ar[d]^{\mu} \\
    \fk^\perp \cap \FN\ar[r] & \FN}  
\end{equation}
where the horizontal arrows are inclusions. Moreover, for $w$ in $W$,
$\dim T^*_w\CB = \dim \CB$ and $T^*_w\CB$ is locally closed in
$Y_{\fk^\perp}$. Thus, the set of irreducible components of
$Y_{\fk^\perp}$ is $\{\, \overline{T^*_w\CB} \mid w\in W\,\}$.

A \emph{quasi-coherent $(\CD_\CB, K)$-module} is a $K$-equivariant,
quasi-coherent $\CD_\CB$-module (for the precise definition see
\cite[\S2]{borhobrylinski:differentialIII}).  If $\CM$ is a coherent
$(\CD_\CB, K)$-module, then $V_{\FNt}(\CM) \subseteq Y_{\fk^\perp}$.

Similarly, a \emph{$(\fg, K)$-module} is a $\fg$-module with a
compatible algebraic action of $K$ (for the precise definition see
\cite[\S2]{borhobrylinski:differentialIII}). If $M$ is a finitely
generated $(\fg, K)$-module, then $V_{\fg}(M)$ is contained in
$\fk^\perp$.

As in the non-equivariant setting, Beilinson-Bernstein
\cite[\S2]{beilinsonbernstein:localisation} have proved that the
global section functor, $\Gamma(\CB, \,\cdot\,)$, defines an
equivalence of categories between the category of quasi-coherent
$(\CD_\CB, K)$-modules and the category of $(\fg, K)$-modules with
trivial central character. Under this equivalence, coherent $(\CD_\CB,
K)$-modules correspond to finitely generated $(\fg, K)$-modules with
trivial central character.

The addition of a $K$-action results in a finer version of Theorem
\ref{thmsigmaa} (see \cite[\S 4]{borhobrylinski:differentialIII}).

\begin{theorem}\label{thmsigmab}
  Suppose that $\CM$ is a coherent $(\CD_\CB, K)$-module and let $M$
  denote the space of global sections of $\CM$.
  \begin{itemize}
  \item[(a)] The variety $V_{\FNt}(\CM)$ is a union of irreducible
    components of $Y_{\fk^\perp}$ and so there is a subset
    $\Sigma(\CM)$ of $W$ such that $V_{\FNt}(\CM)= \bigcup_{w\in
      \Sigma(\CM)} \overline{T^*_w \CB}$.
  \item[(b)] The variety $V_{\fg}(M)$ is contained in $\fk^\perp \cap
    \FN$ and
    \[
    V_{\fg}(M)= \mu(V_{\FNt}(\CM))= \bigcup_{w\in \Sigma(\CM)}
    \mu\left( \overline{T^*_w \CB} \right).
    \]
  \end{itemize}
\end{theorem}

Now it is time to unravel the notation in the highest weight and
Harish-Chandra cases.

First consider the highest weight case when $K=B$. We have $\fk^\perp=
\fb^\perp = \fu$.  Hence, $Y_{\fu^\perp}= \mu\inverse(\fu)\cong \{(x,
B')\in \FN\times \CB \mid x\in \fu\cap \fb'\,\}$. We denote
$Y_{\fu^\perp}$ simply by $Y$ and call it the \emph{conormal variety.}
For $w$ in $W$, $X_w$ is the set of $B$-conjugates of $wBw\inverse$
and $T^*_w \CB\cong \{(x, B')\in \FN\times \CB \mid B'\in X_w,\, x\in
\fu\cap \fb'\,\}$. The projection of $T^*_w \CB$ to $\CB$ is a
$B$-equivariant surjection onto $X_w$ and so $T^*_w \CB \cong
B\times^{B_w} \fu_w$. The diagram (\ref{diagY}) becomes
\[
\xymatrix{Y \ar[r] \ar[d]_{\mu} &\FNt \ar[d]^{\mu} \\
  \fu \ar[r] & \FN.}
\]
Moreover, for $w$ in $W$, $\mu \left( T^*_w \CB \right)= B\fu_w$.
Since $\mu$ is proper, it follows that $\mu \left( \overline{T^*_w
    \CB} \right)= \overline{B\fu_w}$ is the closure of an orbital
variety.

Arguments in the spirit of those given in \S\ref{s2.1} (see \cite[\S3]
{hinichjoseph:orbital}) show that if we set $Y_w= T^*_w \CB$ and
$Y_\FO= \mu\inverse(\FO\cap \fu)$, then $\dim Y_\FO=n$,
$\overline{Y_\FO}$ is equidimensional, and the set of irreducible
components of $\overline{ Y_\FO}$ is $\{\, \overline{Y_\FO \cap Y_w}
\mid w\in W_\FO\,\}$.

Next consider the Harish-Chandra case. In this setting, the ambient
group is $G\times G$ and $K=G_d$ is the diagonal subgroup. Clearly,
$\fk^\perp= \fg_d^\perp= \{\, (x, -x)\mid x\in \fg\,\}$ is isomorphic
to $\fg$ and so
\[
Y_{\fg_d^\perp}= (\mu\times \mu) \inverse (\fg_d^\perp)= \{\, (x,-x,
B', B'') \in \fg\times \fg \times \CB \times \CB \mid x\in \fb' \cap
\fb'' \cap \CN \,\}.
\]
Thus, in this case, $Y_{\fg_d^\perp}$ is clearly isomorphic to the
Steinberg variety and we may identify the restriction of $\mu\times
\mu$ to $Y_{\fg_d^\perp}$ with $\mu_z\colon Z\to \FN$. 
The diagram (\ref{diagY}) becomes
\[
\xymatrix{Z \ar[r] \ar[d]_{\mu_z} &\FNt \times \FNt \ar[d]^{\mu \times
    \mu} \\ \FN \ar[r] & \FN \times \FN}
\]
where the bottom horizontal map is given by $x\mapsto (x,-x)$.
Moreover, for $w$ in $W$,
\[
T^*_w(\CB\times \CB)= \{\,(x,-x, B', B'') \mid (B', B'')\in G(B,
wBw\inverse),\, x\in \fb' \cap \fb'' \cap \CN\,\} \cong Z_w.
\]

Let $p_3\colon Z\to \CB$ be the projection on the third factor. Then
$p_3$ is $G$-equivariant, $G$ acts transitively on $\CB$, and the
fibre over $B$ is isomorphic to $Y$. This gives yet another
description of the Steinberg variety: $Z\cong G\times^B Y$.

Now consider the following three categories:
\begin{itemize}
\item coherent $(\CD_{\CB\times \CB}, G_d)$-modules, $\Mod\,(
  \CD_{\CB\times \CB}, G_d)^\coh$;
\item finitely generated $(\fg\times \fg, G_d)$-modules with trivial
  central character, $\Mod\,(\fg\times \fg, G_d)^\fingen_{0,0}$; and
\item finitely generated $(\fg, B)$-modules with trivial central
  character, $\Mod\,(\fg,B)^\fingen_0$.
\end{itemize}

We have seen that the global section functor defines an equivalence of
categories between $\Mod\,( \CD_{\CB\times \CB}, G_d)^\coh$ and
$\Mod\,(\fg\times \fg, G_d)^\fingen_{0,0}$. Bernstein and Gelfand
\cite{bernsteingelfand:tensor}, as well as Joseph
\cite{joseph:dixmier}, have constructed an equivalence of categories
between $\Mod\,(\fg\times \fg, G_d)^\fingen_{0,0}$ and
$\Mod\,(\fg,B)^\fingen_0$.

Composing these two equivalences of categories we see that the
category of coherent $(\CD_{\CB\times \CB}, G_d)$-modules is
equivalent to the category of finitely generated $(\fg, B)$-modules
with trivial central character, $\Mod\,(\fg,B)^\fingen_0$. Both
equivalences behave well with respect to characteristic varieties and
associated varieties and hence so does their composition. This is the
content of the next theorem. The theorem extends Theorem
\ref{thmsigmab} and summarizes the relationships between the various
constructions in this subsection. See
\cite[\S4]{borhobrylinski:differentialIII} for the proof.

\begin{theorem}
  Suppose $\CM$ is a coherent $(\CD_{\CB\times \CB}, G_d)$-module, $M$
  is the space of global sections of $\CM$, and $L$ is the finitely
  generated $(\fg, B)$-module with trivial central character
  corresponding to $M$. Let $\Sigma=\Sigma(\CM)$ be as in Theorem
  \ref{thmsigmab}. Then when $\mu\times \mu\colon Y_{\fg_d^\perp} \to
  \fg_d^\perp$ is identified with $\mu_z\colon Z\to \FN$ we have:
  \begin{itemize}
  \item[(a)] The characteristic variety of $\CM$ is $V_{T^*(\CB\times
      \CB)} (\CM)= \cup_{y\in \Sigma} \overline{Z_y}$, a union of
    irreducible components of the Steinberg variety.
  \item[(b)] The associated variety of $M$ is $V_{\fg}(M)= \mu_z\left(
      V_{\fg} (\CM)\right) = \cup_{y\in \Sigma} \overline{G\fu_y} =
    G\cdot V_{\fu}(L)$, so the associated variety of $M$ is the image
    under $\mu_z$ of the characteristic variety of $\CM$ and is also
    the $G$-saturation of the associated variety of $L$.
  \item[(c)] The associated variety of $L$ is $V_{\fu}(L)= \cup_{y\in
      \Sigma} \overline{B\fu_y}$, a union of closures of orbital
    varieties.
\end{itemize}
\end{theorem}

The characteristic variety of a coherent $(\CD_{\CB\times \CB},
G_d)$-module is the union of the characteristic varieties of its
composition factors. Similarly the associated variety of a finitely
generated $(\fg\times \fg, G_d)$-module or a finitely generated $(\fg,
B)$-module depends only on its composition factors. Thus, computing
characteristic and associated varieties reduces to the case of simple
modules. The simple objects in each of these categories are indexed by
$W$, see \cite[\S3]{beilinsonbernstein:localisation} and \cite[\S2.7,
4.3, 4.8]{borhobrylinski:differentialIII}. If $w$ is in $W_\FO$ and
$\CM_w$, $M_w$, and $L_w$ are corresponding simple modules, then it is
shown in \cite[\S4.9] {borhobrylinski:differentialIII} that
$\mu_z\left(V_{\fg} (\CM_w)\right) = V(M_w)= G\cdot V(L_w)=
\overline{\FO}$.

In general, explicitly computing the subset $\Sigma= \Sigma(\CM_w)$ so
that $V_Z(\CM_w)= \cup_{y\in \Sigma} \overline{Z_y}$ and $V_\fu(L_w)=
\cup_{y\in \Sigma} \overline{B\fu_y}$ for $w$ in $W$ is a very
difficult and open problem. See
\cite[\S4.3]{borhobrylinski:differentialIII} and
\cite[\S6]{hinichjoseph:orbital} for examples and more information.

\subsection{Generalized Steinberg varieties}\label{s2.4}

When analyzing the restriction of a Springer representation to
parabolic subgroups of $W$, Springer introduced a generalization of
$\FNt$ depending on a parabolic subgroup $P$ and a nilpotent orbit in
a Levi subgroup of $P$. Springer's ideas extend naturally to what we
call ``generalized Steinberg varieties.'' The results in this
subsection may be found in \cite{douglassroehrle:geometry}.

Suppose $\CP$ is a conjugacy class of parabolic subgroups of $G$. The
unipotent radical of a subgroup, $P$, in $\CP$ will be denoted by
$U_P$. A $G$-equivariant function, $c$, from $\CP$ to the power set of
$\FN$ with the properties 
\begin{itemize}
\item[(1)] $\fu_P \subseteq c(P)\subseteq \FN\cap \fp$ and
\item[(2)] the image of $c(P)$ in $\fp/ \fu_P$ is the closure of a
  single nilpotent adjoint $P/U_P$-orbit
\end{itemize}
is called a \emph{Levi class function} on $\CP$. Define 
\[
\FNt^{\CP}_c= \{\,(x, P)\in \FN\times \CP\mid x\in c(P) \,\}.
\]
Let $\mu^\CP_c\colon \FNt^\CP_c\to \FN$ denote the projection on the
first factor. Notice that $\mu^\CP_c$ is a proper morphism.

If $\CQ$ is another conjugacy class of parabolic subgroups of $G$ and
$d$ is a Levi class function on $\CQ$, then the \emph{generalized
  Steinberg variety} determined by $\CP$, $\CQ$, $c$, and $d$ is
\[
\XPQcd= \{\, (x, P,Q)\in \FN\times \CP \times \CQ\mid x\in c(P) \cap
d(Q) \,\} \cong \FNt^{\CP}_c \times_{\FN} \FNt^{\CQ}_d.
\]
Since $G$ acts on $\FN$, $\CP$, and $\CQ$, there is a diagonal action
of $G$ on $\XPQcd$ for all $\CP$, $\CQ$, $c$, and $d$.

The varieties arising from this construction for some particular
choices of $\CP$, $\CQ$, $c$, and $d$ are worth noting.
\begin{itemize}
\item[(1)] When $\CP= \CQ=\CB$, then $c(B')= d(B')=\{ \fu_{B'}\}$ for
  every $B'$ in $\CB$, and so $X^{\CB, \CB}_{0,0}= Z$ is \emph{the}
  Steinberg variety of $G$.
\item[(2)] In the special case when $c(P)$ and $d(Q)$ are as small as
  possible and correspond to the zero orbits in $\fp/\fu_P$ and
  $\fq/\fu_Q$ respectively: $c(P)= \fu_P$ and $d(Q)= \fu_Q$, we denote
  $\XPQcd$ by $\XPQoo$. We have $\XPQoo\cong T^*\CP \times_\FN
  T^*\CQ$.
\item[(3)] When $\CP= \CQ= \{G\}$, $\FO_1$ and $\FO_2$ are two
  nilpotent orbits in $\fg$, $c(G)= \overline{ \FO_1}$ and $d(G)=
  \overline{ \FO_2}$, then $X^{\{G\}, \{G\}} _{c,d}\cong \overline{
    \FO_1}\cap \overline{ \FO_2}$.
\end{itemize}

A special case that will arise frequently in the sequel is when $c(P)$
and $d(Q)$ are as large as possible and correspond to the regular,
nilpotent orbits in $\fp/\fu_P$ and $\fq/\fu_Q$ respectively: $c(P)=
\FN\cap \fp$ and $d(Q)= \FN \cap \fq$. We denote this generalized
Steinberg variety simply by $\XPQ$.

Abusing notation slightly, we let $\mu\colon \XPQcd\to \FN$ denote the
projection on the first coordinate and $\pi\colon \XPQcd\to \CP\times
\CQ$ the projection on the second and third coordinates. We can then
investigate the varieties $\XPQcd$ using preimages of $G$-orbits in
$\FN$ and $\CP\times \CQ$ under $\mu$ and $\pi$ as we did in
\S\ref{s2.1} for the Steinberg variety.  Special cases when at least
one of $c(P)$ or $d(Q)$ is smooth turn out to be the most tractable.
We will describe these cases in more detail below and refer the reader
to \cite{douglassroehrle:geometry} for more general results for
arbitrary $\CP$, $\CQ$, $c$, and $d$.

Fix $P$ in $\CP$ and $Q$ in $\CQ$ with $B\subseteq P\cap Q$. Let $W_P$
and $W_Q$ denote  the Weyl groups of $(P,T)$ and $(Q,T)$ respectively.
We consider $W_P$ and $W_Q$ as subgroups of $W$.

For $B'$ in $\CB$, define $\pi_{\CP} (B')$ to be the unique subgroup
in $\CP$ containing $B'$. Then $\pi_{\CP}\colon \CB\to \CP$ is a proper
morphism with fibres isomorphic to $P/B$. Define
\[
\eta\colon Z\to \XPQ\quad \text{by} \quad \eta(x, B', B'')
= (x, \pi_{\CP}(B'), \pi_{\CQ}(B'')).
\]
Then $\eta$ depends on $\CP$ and $\CQ$ and is a proper,
$G$-equivariant, surjective morphism.

Next, set $\ZPQ=\eta \inverse \left (\XPQoo \right)$ and denote the
restriction of $\eta$ to $\ZPQ$ by $\eta_1$. Then $\eta_1$ is also a
proper, surjective, $G$-equivariant morphism. Moreover, the fibres of
$\eta_1$ are all isomorphic to the smooth, complete variety $P/B
\times Q/B$. More generally, define $\ZPQcd= \eta \inverse \left
  (\XPQcd \right)$.

The various varieties and morphisms we have defined fit together in a
commutative diagram where the horizontal arrows are closed embeddings,
the vertical arrows are proper maps, and the squares are cartesian:
\[
\xymatrix{
\ZPQ \ar[r] \ar[d]_{\eta_1} & \ZPQcd \ar[r] \ar[d]_{\eta} &Z
\ar[d]_{\eta} \\ \XPQoo \ar[r] & \XPQcd \ar[r]& \XPQ.}
\]

For $w$ in $W$, define $\ZPQw$ to be the intersection $\ZPQ \cap Z_w$.
Since $(0,B, wBw\inverse)$ is in $\ZPQw$ and $\eta_1$ is
$G$-equivariant, it is straightforward to check that $\ZPQw \cong
G\times^{B_w} \left( \fu_P \cap w\fu_Q \right)$.  Thus $\ZPQw$ is
smooth and irreducible.

The following statements are proved in
\cite{douglassroehrle:geometry}.
\begin{itemize}
\item[(1)] For $w$ in $W$, $\dim \eta(Z_w)\leq 2n$ with equality if
  and only if $w$ has minimal length in $W_Pw W_Q$. The variety $\XPQ$
  is equidimensional with dimension equal to $2n$ and the set of
  irreducible components of $\XPQ$ is
  \[
  \{\, \eta(\Zwbar) \mid \text{$w$ has minimal length in $W_Pw W_Q$}
  \,\}.
  \]
\item[(2)] For $w$ in $W$, $\ZPQw= Z_w$ if and only if $w$ has maximal
  length in $W_Pw W_Q$. The variety $\ZPQ$ is equidimensional with
  dimension equal to $2n$ and the set of irreducible components of
  $\ZPQ$ is
  \[
  \{\, \Zwbar \mid \text{$w$ has maximal length in $W_Pw W_Q$} \,\}.
  \]
\item[(3)] The variety $\XPQoo$ is equidimensional with dimension
  equal to $\dim \fu_P+ \dim \fu_Q$ and the set of irreducible
  components of $\XPQoo$ is
  \[
  \{\, \eta_1(\Zwbar) \mid \text{$w$ has maximal length in $W_Pw W_Q$}
  \,\}.
  \] 
\item[(4)] For a Levi class function $d$ on $\CQ$, define $\rho_d$ to
  be the number of irreducible components of $d(Q)\cap \left(\fu \cap
    \fl_Q \right)$, where $L_Q$ is the Levi factor of $Q$ that
  contains $T$. Then $\rho_d$ is the number of orbital varieties for
  the open dense $L_Q$-orbit in $d(Q)/\fu_Q$ in the variety of
  nilpotent elements in $\fq/\fu_Q\cong \fl_Q$. The varieties $X^{\CB,
    \CQ}_{0,d}$ are equidimensional with dimension $\frac 12 (\dim \fu
  +\dim d(Q) +\dim \fu_Q)$ and $|W:W_\CQ| \rho_d$ irreducible
  components.
\end{itemize}
Notice that the first statement relates minimal double coset
representatives to regular orbits in Levi subalgebras and the third
statement relates maximal double coset representatives to the zero
orbits in Levi subalgebras.

The quantity $\rho_d$ in the fourth statement is the degree of an
irreducible representation of $W_Q$ (see \S\ref{s3.5}) and so
$|W:W_\CQ| \rho_d$ is the degree of an induced representation of $W$.
The fact that $X^{\CB, \CQ}_{0,d}$ has $|W:W_\CQ| \rho_d$ irreducible
components is numerical evidence for Conjecture \ref{conj} below.

\section{Homology}

In this section we take up the rational Borel-Moore homology of the
Steinberg variety and generalized Steinberg varieties. As mentioned in
the Introduction, soon after Steinberg's original paper, Kazhdan and
Lusztig \cite{kazhdanlusztig:topological} defined an action of
$W\times W$ on the top Borel-Moore homology group of $Z$. They
constructed this action by defining an action of the simple
reflections in $W\times W$ on $H_i(Z)$ and showing that the defining
relations of $W\times W$ are satisfied. They then proved that the
representation of $W\times W$ on $H_{4n}(Z)$ is equivalent to the
two-sided regular representation of $W$, and following a suggestion of
Springer, they gave a decomposition of $H_{4n}(Z)$ in terms of
Springer representations of $W$. Springer representations of $W$ will
be described in \S\ref{s3.4}-- \S\ref{s3.6}.

In the mid 1990s Ginzburg \cite[Chapter
3]{chrissginzburg:representation} popularized a quite general
convolution product construction that defines a $\BBQ$-algebra
structure on $H_*(Z)$, the total Borel-Moore homology of $Z$, 
and a ring structure $K^{\Gbar}(Z)$ (see the
next section for $K^{\Gbar}(Z)$).  With this multiplication,
$H_{4n}(Z)$ is a subalgebra isomorphic to the group algebra of $W$.

In this section, following \cite[Chapter
3]{chrissginzburg:representation}, \cite{douglassroehrle:coinvariant},
and \cite{hinichjoseph:orbital} we will first describe the algebra
structure of $H_*(Z)$, the decomposition of $H_{4n}(Z)$ in terms of
Springer representations, and the $H_{4n}(Z)$-module structure on
$H_{2n}(Y)$ using elementary topological constructions. Then we will
use a more sophisticated sheaf-theoretic approach to give an alternate
description of $H_*(Z)$, a different version of the decomposition of
$H_{4n}(Z)$ in terms of Springer representations, and to describe the
Borel-Moore homology of some generalized Steinberg varieties.

\subsection{Borel-Moore homology and convolution}\label{s3.1} 

We begin with a brief review of Borel-Moore homology, including the
convolution and specialization constructions. The definitions and
constructions in this subsection make sense in a very general setting,
however for simplicity we will consider only complex algebraic
varieties. More details and proofs may be found in \cite[Chapter 2]
{chrissginzburg:representation}.

Suppose that $X$ is a $d$-dimensional, quasi-projective, complex
algebraic variety (not necessarily irreducible). Topological notions
will refer to the Euclidean topology on $X$ unless otherwise
specified. Two exceptions to this convention are that we continue to
denote the dimension of $X$ as a complex variety by $\dim X$ and that
``irreducible'' means irreducible with respect to the Zariski
topology. In particular, the topological dimension of $X$ is $2\dim
X$.

Let $X\cup \{\infty\}$ be the one-point compactification of $X$. Then
the $i\th$ Borel-Moore homology space of $X$, denoted by $H_i(X)$, is
defined by $H_i(X)=H_i^{\operatorname{sing}}(X\cup \{ \infty \},
\{\infty\})$, the relative, singular homology with rational
coefficients of the pair $(X\cup \{\infty\}, \{\infty\})$. Define a
graded $\BBQ$-vector space, 
\[
H_*(X)= \sum_{i\geq 0} H_i(X)\ \text{-- the \emph{Borel-Moore homology
    of $X$.}}
\]

Borel-Moore homology is a bivariant theory in the sense of Fulton and
MacPherson \cite{fultonmacpherson:categorical}: Suppose that
$\phi\colon X\to Y$ is a morphism of varieties. 
\begin{itemize}
\item If $\phi$ is proper, then there is an induced direct image map in
  Borel-Moore homology, $\phi_*:H_i(X)\to H_i(Y)$.
\item If $\phi$ is smooth with $f$-dimensional fibres, then there is a
  pullback map in Borel-Moore homology, $\phi^*\colon H_i(Y)\to
  H_{i+2f}(X)$.
\end{itemize}

Moreover, if $X$ is smooth and $A$ and $B$ are closed subvarieties of
$X$, then there is an intersection pairing $\cap\colon H_i(A) \times
H_j(B)\to H_{i+j-2d}(A\cap B)$. Although not reflected in the
notation, this pairing depends on the triple $(X,A,B)$. In particular,
the intersection pairing depends on the smooth ambient variety $X$.
 
In dimensions greater than or equal $2\dim X$, the Borel-Moore
homology spaces of $X$ are easily described. If $i>2d$, then
$H_i(X)=0$, while the space $H_{2d}(X)$ has a natural basis indexed by
the $d$-dimensional irreducible components of $X$. If $C$ is a
$d$-dimensional irreducible component of $X$, then the homology class
in $H_{2d}(X)$ determined by $C$ is denoted by $[C]$.

For example, for the Steinberg variety, $H_i(Z)=0$ for $i>4n$ and the
set $\{\, [\Zwbar]\mid w\in W\,\}$ is a basis of $H_{4n}(Z)$.
Similarly, for the conormal variety, $H_i(Y)=0$ for $i>2n$ and the set
$\{\, [\Ywbar]\mid w\in W\,\}$ is a basis of $H_{2n}(Y)$.

Suppose that for $i=1,2,3$, $M_i$ is a smooth, connected,
$d_i$-dimensional variety. For $1\leq i<j\leq 3$, let $p_{i,j}\colon
M_1\times M_2 \times M_3\to M_i\times M_j$ denote the projection.
Notice that each $p_{i,j}$ is smooth and so the pullback maps
$p_{i,j}^*$ in Borel-Moore homology are defined.

Now suppose $Z_{1,2}$ is a closed subset of $M_1\times M_2$ and $Z_{2,3}$
is a closed subvariety of $M_2\times M_3$. Define $Z_{1,3} = Z_{1,2}
\circ Z_{2,3}$ to be the composition of the relations $Z_{1,2}$ and
$Z_{2,3}$. Then
\[
Z_{1,3} =\{\, (m_1, m_3)\in M_1\times M_3\mid \text{$\exists\, m_2\in
  M_2$ with $(m_1, m_2)\in Z_{1,2}$ and $(m_2, m_3)\in Z_{2,3}$} \,\}.
\]

In order to define the convolution product, we assume in addition that
the restriction
\[
p_{1,3}\colon p_{1,2}\inverse(Z_{1,2}) \cap p_{2,3}
\inverse (Z_{2,3}) \to Z_{1,3}
\]
is a proper morphism. Thus, there is a direct image map
\[
(p_{1,3})_*\colon H_i \left( p_{1,2} \inverse(Z_{1,2}) \cap p_{2,3}
  \inverse (Z_{2,3}) \right) \to H_{i} (Z_{1,3})
\]
in Borel-Moore homology. The \emph{convolution product}, $H_i(Z_{1,2})
\times H_j(Z_{2,3})\xrightarrow{*} H_{i+j-2d_2} (Z_{1,3})$ is then
defined by
\[
c*d= (p_{1,3})_*\left( p_{1,2}^*(c) \cap p_{2,3}^*(d) \right)
\]
where $\cap$ is the intersection pairing determined by the subsets
$Z_{1,2} \times M_3$ and $M_1 \times Z_{2,3}$ of $M_1\times M_2\times
M_3$. It is a straightforward exercise to show that the convolution
product is associative.

The convolution construction is particularly well adapted to fibred
products. Fix a ``base'' variety, $N$, which is not necessarily
smooth, and suppose that for $i=1,2,3$, $f_i\colon M_i\to N$ is a
proper morphism. Then taking $Z_{1,2}= M_1\times_N M_2$, $Z_{2,3}=
M_2\times_N M_3$, and $Z_{1,3}= M_1\times_N M_3$, we have a
convolution product $H_i(M_1\times_N M_2) \times H_j(M_2\times_N
M_3)\xrightarrow{*} H_{i+j-2d_2} (M_1\times_N M_3)$.

As a special case, when $M_1=M_2=M_3=M$ and $f_1=f_2=f_3=f$, then
taking $Z_{i,j}= M\times_N M$ for $1\leq i<j\leq 3$, the convolution
product defines a multiplication on $H_*(M\times_N M)$ so that
$H_*(M\times_N M)$ is a $\BBQ$-algebra with identity. The identity in
$H_*(M\times_N M)$ is $[M_{\Delta}]$ where $M_{\Delta}$ is the
diagonal in $M\times M$. If $d=\dim M$, then $H_i(M\times_N M)*
H_j(M\times_N M)\subseteq H_{i+j-2d}(M\times_N M)$ and so
$H_{2d}(M\times_N M)$ is a subalgebra and $\oplus_{i<2d} H_i(M\times_N
M)$ is a nilpotent, two-sided ideal.

Another special case is when $M$ and $M'$ are smooth and $f\colon M\to
N$ and $f'\colon M'\to N$ are proper maps. Then taking $Z_{1,2}=
M\times_NM$ and $Z_{2,3}= M\times_NM'$, the convolution product
defines a left $H_*(M\times_N M)$-module structure on $H_*(M\times_N
M')$. A further special case of this construction is when $M'=A$ is a
smooth, closed subset of $N$ and $f'\colon A\to N$ is the inclusion.
Then $M\times_NA \cong f\inverse A$ and the convolution product
defines a left $H_*(M\times_N M)$-module structure on
$H_*(f\inverse(A))$. This construction will be exploited extensively
in \S\ref{s3.5}.

As an example, recall that $Z\cong \FNt\times_\FN \FNt$ where
$\mu\colon \FNt\to \FN$ is a proper map. Applying the constructions in
the last two paragraphs to $Z$ and to $M'$, where
$M'=Y=\mu\inverse(\fu)$ and $M'=\CB_x=\mu\inverse(x)$ for $x$ in
$\FN$, we obtain the following proposition.

\begin{proposition}
  The convolution product defines a $\BBQ$-algebra structure on
  $H_*(Z)$ so that $H_{4n}(Z)$ is a $|W|$-dimensional subalgebra and
  $\bigoplus_{i<4n} H_i(Z)$ is a two-sided, nilpotent ideal. Moreover,
  the convolution product defines left $H_*(Z)$-module structures on
  $H_*(Y)$ and on $H_*(\CB_x)$ for $x$ in $\FN$.
\end{proposition}

In the next two subsections we will make use of the following
specialization construction in Borel-Moore homology due to Fulton and
MacPherson \cite[\S 3.4]{fultonmacpherson:categorical}.

Suppose that our base variety $N$ is smooth and $s$-dimensional. Fix a
distinguished point $n_0$ in $N$ and set $N^*=N\setminus \{n_0\}$.
Let $M$ be a variety, not necessarily smooth, and suppose that
$\phi\colon M\to N$ is a surjective morphism. Set $M_0=
\phi\inverse(n_0)$ and $M^*=\phi \inverse(N^*)$. Assume that the
restriction $\phi|_{M^*} \colon M^* \to N^*$ is a locally trivial
fibration. Then there is a ``specialization'' map in Borel-Moore
homology, $\lim\colon H_i(M^*)\to H_{i-2s}(M_0)$ (see
\cite[\S2.6]{chrissginzburg:representation}). It is shown in
\cite[\S2.7]{chrissginzburg:representation} that when all the various
constructions are defined, specialization commutes with convolution:
$\lim (c*d)= \lim c *\lim d$.

\subsection{The specialization construction and $H_{4n}(Z)$}
\label{s3.2} 

Chriss and Ginzburg \cite[\S3.4]{chrissginzburg:representation} use
the specialization construction to show that $H_{4n}(Z)$ is isomorphic
to the group algebra $\BBQ [W]$. We present their construction in this
subsection. In the next subsection we show that the specialization
construction can also be used to show that $H_*(Z)$ is isomorphic to
the smash product of the group algebra of $W$ and the coinvariant
algebra of $W$.

We would like to apply the specialization construction when the
variety $M_0$ is equal $Z$. In order to do this, we need varieties
that are larger than $\FN$, $\FNt$, and $Z$.

Define 
\[
\fgt= \{\,(x, B')\in \fg\times \CB\mid x\in \fb'\,\} \quad \text{and}
\quad \Zhat=\{\, (x, B', B'')\in \fg\times \CB \times \CB \mid x\in
\fb' \cap \fb''\,\}.
\]
Abusing notation again, let $\mu\colon \fgt\to \fg$ and $\mu_z\colon
\Zhat\to \fg$ denote the projections on the first factors and let
$\pi\colon \Zhat\to \CB \times \CB$ denote the projection on the
second and third factors.

For $w$ in $W$ define $\Zhat_w= \pi\inverse(G(B, wBw\inverse))$. Then
$\Zhat_w\cong G\times^{B_w} \fb_w$. Therefore, $\dim \Zhat_w= \dim
\fg$ and the closures of the $\Zhat_w$'s for $w$ in $W$ are the
irreducible components of $\Zhat$.

As with $Z$, we have an alternate description of $\Zhat$ as $(\fgt
\times \fgt) \times_{\fg\times \fg} \fg$. However, in contrast to the
situation in \S\ref{s2.3}, where $Z\cong \{\, (x,-x, B',B'')\in \FN
\times \FN\times \CB \times \CB \mid x\in \fb' \cap \fb'' \cap \FN
\,\}$, in this section we use that $\Zhat \cong \{\, (x,B', x,B'')\in
\fg \times \CB\times \fg \times \CB \mid x\in \fb' \cap \fb''\,\}$. In
particular, we will frequently identify $\Zhat$ with the subvariety of
$\fgt\times \fgt$ consisting of all pairs $((x,B'), (x,B''))$ with $x$
in $\fb'\cap \fb''$. Similarly, we will frequently identify $Z$ with
the subvariety of $\FNt\times \FNt$ consisting of all pairs $((x,B'),
(x,B''))$ with $x$ in $\FN \cap \fb' \cap \fb''$.

For $(x, gBg\inverse)$ in $\fgt$, define $\nu (x,gBg\inverse)$ to be
the projection of $g\inverse\cdot x$ in $\ft$. Then $\nu\colon \fgt\to
\ft$ is a surjective morphism. For $w$ in $W$, let
$\Gamma_{w\inverse}= \{\, (h, w\inverse \cdot h)\mid h\in \ft\,\}
\subseteq \ft \times \ft$ denote the graph of the action of
$w\inverse$ on $\ft$ and define
\[
\Lambda_w= \Zhat \cap (\nu\times \nu)\inverse \left(
  \Gamma_{w\inverse} \right) =\{\, (x, B', B'')\in \Zhat \mid
\nu(x,B'')= w\inverse \nu(x,B')\,\}.
\]

The spaces we have defined so far fit into a commutative diagram with
cartesian squares where $\delta\colon \fg \to \fg \times \fg$ is the
diagonal map:
\begin{equation}
  \label{eq:nu}
  \xymatrix{
    \Lambda_w \ar[r] \ar[d] & \Zhat \ar[r]^{\mu_z} \ar[d] & \fg
    \ar[d]^{\delta} \\ 
    (\nu\times \nu)\inverse \left(\Gamma_{w\inverse} \right) \ar[r]
    \ar[d]  & \fgt \times \fgt \ar[r]_{\mu\times \mu}
    \ar[d]_{\nu\times \nu} & \fg \times \fg \\ 
    \Gamma_{w\inverse} \ar[r] & \ft\times \ft.&}  
\end{equation}

Let $\nu_w\colon \Lambda_w\to \Gamma_{w\inverse}$ denote the
composition of the leftmost vertical maps in (\ref{eq:nu}), so $\nu_w$
is the restriction of $\nu\times \nu$ to $\Lambda_w$.  We will
consider subsets of $\Zhat$ of the form $\nu_w\inverse(S')$ for
$S'\subseteq \Gamma_{w\inverse}$.  Thus, for $h$ in $\ft$ we define
$\Lambda_w^h= \nu_w\inverse(h, w\inverse h)$. Notice in particular
that $\Lambda_w^0= Z$. More generally, for a subset $S$ of $\ft$ we
define $\Lambda_w^{S}=\coprod _{h\in S} \Lambda_w^h$. Then
$\Lambda_w^{S}= \nu_w\inverse(S')$, where $S'$ is the graph of
$w\inverse$ restricted to $S$.

Let $\ft_\reg$ denote the set of regular elements in $\ft$. For $w$
in $W$, define $\wtilde\colon G/T\times \ft_\reg\to G/T\times
\ft_\reg$ by $\wtilde(gT, h)= (gwT, w\inverse h)$. The rule $(gT,
h)\mapsto (g\cdot h, gB)$ defines an isomorphism of varieties $f\colon
G/T\times \ft_\reg \xrightarrow{\cong} \fgrst$, where $\fgrst=
\mu\inverse(G \cdot \ft_\reg)$. We denote the automorphism $f\circ
\wtilde \circ f\inverse$ of $\fgrst$ also by $\wtilde$.

We now have all the notation in place for the specialization
construction. Fix an element $w$ in $W$ and a one-dimensional
subspace, $\ell$, of $\ft$ so that $\ell\cap \ft_\reg= \ell\setminus
\{0\}$.  The line $\ell$ is our base space and the distinguished point
in $\ell$ is $0$. As above, we set $\ell^*= \ell\setminus \{0\}$. We
denote the restriction of $\nu_w$ to $\Lambda_w^{\ell}$ again by
$\nu_w$. Then $\nu_w\colon \Lambda_w^{\ell}\to \ell$ is a surjective
morphism with $\nu_w\inverse (0)=Z$ and $\nu_w\inverse( \ell^*)
=\Lambda_w^{\ell^*}$. We will see below that the restriction
$\Lambda_w^{\ell^*} \to \ell^*$ is a locally trivial fibration and so
a specialization map
\begin{equation}
  \label{eq:lim}
  \lim\colon H_{i+2}( \Lambda_w^{\ell^*}) \to H_{i}(Z)
\end{equation}
is defined.

It is not hard to check that the variety $\Lambda_w^{\ell^*}$ is the
graph of $\wtilde|_{ \fgt^{\ell^*}}\colon \fgt^{\ell^*}\to
\fgt^{w\inverse(\ell^*)}$, where for an arbitrary subset $S$ of $\ft$,
$\fgt^S$ is defined to be $\nu\inverse(S)= \{\, (x, B')\in \fgt \mid
\nu(x,B')\in S\,\}$. It follows that for $h$ in $\ell^*$ we have
$\nu_w\inverse(h) = \Lambda_w^h \cong G/T$ and that
$\Lambda_w^{\ell^*}\to \ell^*$ is a locally trivial fibration.
Moreover, $\Lambda_w^{\ell^*} \cong \fgt^{\ell^*}$ and hence is an
irreducible, $(2n+1)$-dimensional variety. Therefore, $H_{4n+2}(
\Lambda_w^{ \ell^*})$ is one-dimensional with basis $\{
[\Lambda_w^{\ell^*}] \}$.  Taking $i=4n$ in (\ref{eq:lim}), we define
\[
\lambda_w= \lim ([\Lambda_w^{\ell^*}])
\]
in $H_{4n}(Z)$.

Because $\Lambda_w^{\ell^*}$ is a graph, it follows easily from the
definitions that for $y$ in $W$, there is a convolution product
\[
H_*(\Lambda_{w}^{\ell^*}) \times H_*(\Lambda_{y}^{w\inverse\ell^*})
\xrightarrow{*} H_*(\Lambda_{wy}^{\ell^*})
\]
and that $[\Lambda_{w}^{\ell^*}] *[\Lambda_{y}^{w\inverse\ell^*}] =
[\Lambda_{wy}^{\ell^*}]$. Because specialization commutes with
convolution, we have $\lambda_w * \lambda_y= \lambda_{wy}$ for all $w$
and $y$ in $W$.

Chriss and Ginzburg \cite[\S3.4]{chrissginzburg:representation} have
proved the following:
\begin{itemize}
\item[(1)] The element $\lambda_w$ in $H_{4n}(Z)$ does not depend on
  the choice of $\ell$.
\item[(2)] The expansion of $\lambda_w$ as a linear combination of the
  basis elements $[\overline{Z_y}]$ of $H_{4n}(Z)$ has the form
  $\lambda_w= [\Zwbar] +\sum_{y< w} a_{y,w} [\overline{Z_y}]$ where
  $\leq$ is the Bruhat order on $W$.
\end{itemize}
These results prove the following theorem.

\begin{theorem}\label{thm3.2}
  With the notation as above, the assignment $w\mapsto \lambda_w$
  extends to an algebra isomorphism $\BBQ [W] \xrightarrow{\cong}
  H_{4n}(Z)$.
\end{theorem}

\subsection{The Borel-Moore homology of $Z$ and coinvariants}

Now consider
\[
Z_1 =\{\, (x, B', B')\in \FN\times \CB\times \CB\mid x\in \fb'\,\}.
\]
Then $Z_1$ may be identified with the diagonal in $\FNt\times \FNt$.
It follows that $Z_1$ is closed in $Z$ and isomorphic to $\FNt$.

Since $\FNt \cong T^*\CB$, it follows from the Thom isomorphism in
Borel-Moore homology that $H_{i+2n}(Z_1) \cong H_{i}(\CB)$ for all
$i$. Since $\CB$ is smooth and compact, $H_{i}(\CB)\cong
H^{2n-i}(\CB)$ by Poincar\'e duality. Therefore, $H_{4n-i}(Z_1) \cong
H^i(\CB)$ for all $i$.

The cohomology of $\CB$ is well-understood: there is an isomorphism of
graded algebras, $H^{*}(\CB)\cong \Coinv_*(W)$ where $\Coinv_*(W)$ is
the coinvariant algebra of $W$ with generators in degree 2.  It
follows that $H_j(Z_1)=0$ if $j$ is odd, $H_{4n-2i}(Z_1)\cong
\Coinv_{2i}(W)$ for $0\leq i\leq n$.

The following is proved in \cite{douglassroehrle:coinvariant}.
\begin{itemize}
\item[(1)] There is a convolution product on $H_*(Z_1)$. With this
  product, $H_{*}(Z_1)$ is a commutative $\BBQ$-algebra and there is
  an isomorphism of graded $\BBQ$-algebras
  \[
  \beta\colon \Coinv_*(W) \xrightarrow{\cong} H_{4n-*}
  (Z_1).
  \]
\item[(2)] If $\iota \colon Z_1\to Z$ denotes the inclusion, then the
  direct image map in Borel-Moore homology, $\iota_*\colon H_*(Z_1)
  \to H_*(Z)$, is an injective ring homomorphism.
\item[(3)] If we identify $H_*(Z_1)$ with its image in $H_*(Z)$ as in
  (b), then the linear transformation given by the convolution product
  \[
  H_i(Z_1)\otimes H_{4n}(Z) \xrightarrow{\ *\ } H_i(Z)
  \]
  is an isomorphism of vector spaces for $0\leq i\leq 4n$.
\end{itemize}

The algebra $\Coinv_*(W)$ has a natural action of $W$ by algebra
automorphisms and the isomorphism $\beta$ in (a) is in fact an
isomorphism of $W$-algebras. The $W$-algebra structure on $H_*(Z_1)$
is described as follows.

Fix $w$ in $W$ and identify $H_*(Z_1)$ with its image in $H_*(Z)$.
Then
\begin{equation*}
  \lambda_w*H_i(Z_1) *\lambda_{w\inverse} = H_i(Z_1).
\end{equation*}
Therefore, conjugation by $\lambda_w$ defines a $W$-algebra structure
on $H_*(Z_1)$. With this $W$-algebra structure, the isomorphism
$\beta\colon \Coinv_*(W) \xrightarrow{\cong} H_{4n-*} (Z_1)$ in (a) is
an isomorphism of $W$-algebras.

Using the natural action of $W$ on $\Coinv(W)$, we can define the
smash product algebra $\Coinv(W) \rtimes \BBQ [W]$. We suppose that
$\Coinv(W) \rtimes \BBQ [W]$ is graded by $(\Coinv(W) \rtimes \BBQ
[W])_i= \Coinv_i(W) \otimes \BBQ [W]$. Then combining Theorem
\ref{thm3.2}, item (3) above, and the fact that $\beta$ is an
isomorphism of $W$-algebras, we obtain the following theorem giving an
explicit description of the structure of $H_{*}(Z)$.

\begin{theorem}\label{thm3.3}
  The composition
  \[
  \Coinv_*(W) \rtimes \BBQ [W] \xrightarrow{\beta\otimes \alpha}
  H_{4n-*}(Z_1) \otimes H_{4n}(Z) \xrightarrow{*} H_{4n-*}(Z)
  \]
  is an isomorphism of graded $\BBQ$-algebras.
\end{theorem}

\subsection{Springer representations of $W$} \label{s3.4}

Springer \cite{springer:trigonometric} \cite{springer:construction}
has given a case-free construction of the irreducible representations
of $W$. He achieves this by defining an action of $W$ on $H^*(\CB_x)$
for $x$ in $\FN$.  Define $d_x= \dim \CB_x$ and let $C(x)= Z_G(x)/
Z_G^0(x)$. Then the centralizer in $G$ of $x$ acts on $\CB_x$ and so
$C(x)$ acts on $H^*(\CB_x)$. Springer shows that if $\phi$ is an
irreducible representation of $C(x)$ and $H^{2d_x}(\CB_x)_\phi$ is the
homogeneous component of $H^{2d_x}(\CB_x)$ corresponding to $\phi$,
then $H^{2d_x}(\CB_x)_\phi$ is $W$-stable and is either zero or
affords an irreducible representation of $W$. He shows furthermore
that every irreducible representation of $W$ arises in this way.

We have seen in \S\ref{s3.1} that for $x$ in $\FN$, the convolution
product defines a left $H_{4n}(Z)$-module structure on $H_*(\CB_x)$
and in \S\ref{s3.2} that $H_{4n}(Z) \cong \BBQ[W]$. Thus, we obtain a
representation of $W$ on $H_*(\CB_x)$. Because $\CB_x$ is projective,
and hence compact, $H^*(\CB_x)$ is the linear dual of $H_*(\CB_x)$ and
so we obtain a representation of $W$ on $H^*(\CB_x)$. 

In the next subsection we use topological techniques to decompose
the two-sided regular representation of $H_{4n}(Z)$ into irreducible
sub-bimodules and describe these sub-bimodules explicitly in terms of
the irreducible $H_{4n}(Z)$-submodules of $H_{2d_x}( \CB_x)$ for $x$
in $\FN$.  In \S\ref{s3.6} we use sheaf theoretic techniques to
decompose the representation of $\BBQ[W] \cong H_{4n}(Z)$ on
$H_*(\CB_x)$ into irreducible constituents.

As above, the component group $C(x)$ acts on $H_*(\CB_x)$. It is easy
to check that the $C(x)$-action and the $H_{4n}(Z)$-action commute.
Therefore, up to isomorphism, the representation of $W$ on
$H_{*}(\CB_x)$ depends only on the $G$-orbit of $x$ and the isotopic
components for the $C(x)$-action afford representations of $W$.

It follows from results of Hotta \cite{hotta:springer} that the
representations of $W$ on $H_*(\CB_x)$ constructed using the
convolution product and the isomorphism $\BBQ[W] \cong H_{4n}(Z)$ are
the same as the representations originally constructed by Springer
tensored with the sign representation of $W$.

As an example, consider the special case corresponding to the trivial
representation of $C(x)$: $H_{2d_x}(\CB_x)^{C(x)}$, the
$C(x)$-invariants in $H_{2d_x}(\CB_x)$. Let $\CI rr_x$ denote the set
of irreducible components of $\CB_x$. Then $\{\, [C] \mid C\in \CI
rr_x\,\}$ is a basis of $H_{2d_x} (\CB_x)$. The group $C(x)$ acts on
$H_{2d_x}(\CB_x)$ by permuting this basis: $g [C]=[gC]$ for $g$ in
$Z_G(x)$ and $C$ in $\CI rr_x$. Thus, the orbit sums index a basis of
$H_{2d_x}(\CB_x)^{C(x)}$. We have seen in \S\ref{s2.2} that there is a
bijection between the orbits of $C(x)$ on $\CI rr_x$ and the set of
orbital varieties for $\FO$ where $\FO$ is the $G$-orbit of $x$. Thus,
$H_{2d_x}(\CB_x)^{C(x)}$ affords a representation of $W$ and has a
basis naturally indexed by the set of orbital varieties for $\FO$. It
follows from the general results stated above and discussed in more
detail in the following two subsections that this representation is
irreducible.

\subsection{More on the top Borel-Moore homology of $Z$} \label{s3.5} 

We saw in Theorem \ref{thm3.2} that $H_{4n}(Z) \cong \BBQ[W]$. In this
subsection we follow the argument in
\cite[\S3.5]{chrissginzburg:representation}. First we obtain a
filtration of $H_{4n}(Z)$ by two sided ideals indexed by the set of
nilpotent orbits in $\FN$ and then we describe the decomposition of
the associated graded ring into minimal two-sided ideals. In
particular, we obtain a case-free construction and parametrization of
the irreducible representations of $W$. As explained in the
introduction, a very similar result was first obtained using different
methods by Kazhdan and Lusztig \cite{kazhdanlusztig:topological},
following an idea of Springer.

Recall that orbit closure defines a partial order on the set of
nilpotent orbits in $\FN$: $\FO_1\leq \FO_2$ if $\FO_1\subseteq
\overline{\FO_2}$. For a nilpotent orbit, $\FO$, define $\partial
\FObar= \FObar \setminus \FO= \{\, \FO'\mid \FO'< \FO\,\}$ and set
$Z_{\FObar}= \mu_z\inverse (\FObar)$, and $Z_{\partial \FObar}=
\mu_z\inverse (\partial \FObar)$. Notice that $\partial \FO$ is a
closed subvariety of $\FN$. Define $W_{\FObar} = \cup_{ \FD\subseteq
  \FObar} W_{\FD}$ and $W_{\partial \FObar} = \cup_{ \FD\subseteq
  \partial \FObar} W_{\FD}$, where the union is taken over the
nilpotent orbits contained in $\FObar$ and $\partial \FObar$
respectively.

It follows from the results in \S\ref{s2.1} and \S\ref{s3.1} that
$\{\, [\Zwbar]\mid w\in W_{\FObar} \,\}$ is a basis of
$H_{4n}(Z_{\FObar})$.

If we take $f_i\colon M_i\to N$ to be $\mu\colon \FNt\to \FN$ for
$i=1,2,3$ and $Z_{i,j}= Z_{\FObar}$ for $1\leq i\ne j\leq 3$, then the
convolution product construction in \S\ref{s3.1} defines the structure
of a $\BBQ$-algebra on $H_{*} (Z_{\FObar})$ and $H_{4n} (Z_{\FObar})$
is a subalgebra. Similarly, taking $Z_{1,2}=Z$ and $Z_{2,3}=Z_{1,3}=
Z_{\FObar}$, the convolution product defines a left $H_{*}(Z)$-module
structure on $H_{*}(Z_{\FObar})$ that is compatible with the algebra
structure on $H_{*}(Z_{\FObar})$ in the sense that $a*(y*z)= (a*y)*z$
for $a$ in $H_{*}(Z)$ and $y$ and $z$ in $H_{*}(\FObar)$. Taking
$Z_{1,2}=Z_{1,3}=Z_{\FObar}$ and $Z_{1,2}=Z$, we get a right
$H_{*}(Z)$-module structure on $H_{*}(Z_{\FObar})$ that commutes with
the left $H_{*}(Z)$-module structure and is compatible with the
algebra structure. Thus, we see that $H_{4n}(Z_{\FObar})$ is a
$|W_{\FObar}|$-dimensional algebra with a compatible
$H_{4n}(Z)$-bimodule structure.

Arguing as in the last two paragraphs with $Z_{\FObar}$ replaced by
$Z_{\partial \FObar}$, we see that $H_{4n}(Z_{\partial \FObar})$ is a
$|W_{\partial \FObar}|$-dimensional algebra with a compatible
$H_{4n}(Z)$-bimodule structure.

The inclusions $Z_{\partial \FObar} \subseteq Z_{\FObar} \subseteq Z$
induce injective, $H_{4n}(Z) \times H_{4n}(Z)$-linear ring
homomorphisms, $H_{4n}( Z_{\partial \FObar}) \to H_{4n} (Z_{\FObar})
\to H_{4n}(Z)$, and so we may identify $H_{4n}( Z_{\partial \FObar})$
and $H_{4n}( Z_{\FObar})$ with their images in $H_{4n}(Z)$ and
consider $H_{4n}( Z_{\partial \FObar})$ and $H_{4n}( Z_{\FObar})$ as
two-sided ideals in $H_{4n}(Z)$.

The two-sided ideals $H_{4n}(Z_{\FObar})$ define a filtration of
$H_{4n}(Z)$ indexed by the set of nilpotent orbits.  Thus, to describe
the decomposition of the associated graded algebra into minimal
two-sided ideals, we need to analyze the quotients $H_{4n}(
Z_{\FObar}) /H_{4n}( Z_{\partial \FObar})$. Because $H_{4n}(Z)$ is
semisimple (it is isomorphic to $\BBQ[W]$), this will also describe
the two-sided regular representation of $H_{4n}(Z)$ into minimal
sub-bimodules and give a case-free construction of the irreducible
representations of $W$.

For a $G$-orbit, $\FO$, define $H_{\FO}$ to be the quotient $H_{4n}(
Z_{\FObar}) /H_{4n}( Z_{\partial \FObar})$. Then $\dim H_{\FO} =
|W_\FO|$ and $H_{\FO}$ is an $H_{4n}(Z)$-bimodule with a compatible
$\BBQ$-algebra structure inherited from the convolution product on
$H_{4n}(Z)$.

Now fix a $G$-orbit $\FO$ and an element $x$ is in $\FO$. Set $Z_x=
\mu_z\inverse(x)$. Then clearly $Z_x\cong \CB_x \times \CB_x$ and
$\dim Z_x=2d_x$. The centralizer of $x$ acts diagonally on $Z_x$, and
so the component group, $C(x)$, acts on $H_{*}(Z_x)$. Thus,
$H_{4d_x}(Z_x) ^{C(x)} \cong H_{4d_x}( \CB_x \times \CB_x)^{ C(x)}$
has a basis indexed by the $C(x)$-orbits on the set of irreducible
components of $\CB_x\times \CB_x$. We saw in \S\ref{s2.1} that there
is a bijection between the $C(x)$-orbits on the set of irreducible
components of $\CB_x\times \CB_x$ and the two-sided Steinberg cell
corresponding to $\FO$.  Therefore, the dimension of
$H_{4d_x}(Z_x)^{C(x)}$ is $|W_\FO| =\dim H_{\FO}$.

As for $Z_{\FObar}$ and $Z_{\partial \FObar}$, if we take $f_i\colon
M_i\to N$ to be $\mu\colon \FNt\to \FN$ for $i=1,2,3$, then for
suitable choices of $Z_{i,j}$ for $1\leq i<j\leq 3$, the convolution
product defines a $\BBQ$-algebra structure and a compatible
$H_{*}(Z)$-bimodule structure on $H_{4d_x}(Z_x)$.  It is
straightforward to check that $H_{4d_x}(Z_x)^{C(x)}$ is a subalgebra
and an $H_{*}(Z)$-sub-bimodule of $H_{4d_x}(Z_x)$. 

The group $C(x)$ acts diagonally on $H_{2d_x}(\CB_x) \otimes
H_{2d_x}(\CB_x)$ and it follows from the K\"unneth formula that
\begin{equation}
  \label{kunneth}
  H_{4d_x}(Z_x)^{C(x)} \cong \left( H_{2d_x}(\CB_x) \otimes
    H_{2d_x}(\CB_x) \right) ^{C(x)}.
\end{equation} 
The convolution product defines left and right $H_{*}(Z)$-module
structures on $H_*(\CB_x)$ and the isomorphism in (\ref{kunneth}) is
as $H_{*}(Z)$-bimodules, where $H_{*}(Z)$ acts on the right-hand side
by acting on the left on the first $H_{2d_x}(\CB_x)$ and on the right
on the second $H_{2d_x}(\CB_x)$.

Fix a set, $\FS$, of $G$-orbit representatives in $\FN$. The next
proposition has been proved by Kazhdan and Lusztig
\cite{kazhdanlusztig:topological} and Chriss and Ginzburg
\cite[\S3.5]{chrissginzburg:representation}. An alternate argument has
also been given by Hinich and Joseph \cite[\S4]{hinichjoseph:orbital}.

\begin{proposition}\label{ideal}
  There is an algebra isomorphism $H_{\FO} \cong H_{4n}(Z_x)^{C(x)}$
  and $H_{4n}(Z)$-bimodule isomorphisms
  \[
  H_{\FO} \cong H_{4n}(Z_x)^{C(x)} \cong (H_{2d_x} (\CB_x) \otimes
  H_{2d_x} (\CB_x))^{C(x)}.
  \]
\end{proposition}

For $\FO=\{0\}$, the $H_{4n}(Z)$-bimodule $H_{\FO}$ corresponds to the
trivial representation of $W$ under the isomorphism $H_{4n}(Z) \cong
\BBQ[W]$. For $\FO$ the regular nilpotent orbit, the
$H_{4n}(Z)$-bimodule $H_{\FO}$ corresponds to the sign representation
of $W$. In general however, $H_{\FO}$ is not a minimal two-sided ideal
in the associated graded ring, $\gr\, H_{4n}(Z)$, and not an
irreducible $H_{4n}(Z)$-bimodule. To obtain the decomposition of $\gr
\,H_{4n}(Z)$ into irreducible $H_{4n}(Z)$-bimodules, we need to
decompose each $H_{2d_x}(\CB_x)$ into $C(x)$-isotopic components.

For an irreducible representation of $C(x)$ with character $\phi$,
denote the $\phi$-isotopic component of $C(x)$ on
$H_{2d_x}(\CB_x)$ by $H_{2d_x}(\CB_x)_\phi$. Define $\widehat{C(x)}$
to be the set of $\phi$ with $H_{2d_x}(\CB_x)_\phi\ne 0$. We saw in
the last subsection that the trivial character of $C(x)$ is always an
element of $\widehat{C(x)}$. The sets $\widehat{C(x)}$ have been
computed explicitly in all cases, see \cite[\S13.3]{carter:finite}.
For example, if $G=\GL_n(\BBC)$, then $Z_G(x)$ is connected and so 
$C(x) = 1$ for all $x$ in $\FN$, and so $\widehat{C(x)}$ contains all irreducible
characters of $C(x)$.  In general $\widehat{C(x)}$ does not contain
all irreducible characters of $C(x)$. 

Recall from \S\ref{s3.4} that for each $\phi$, $H_{2d_x}(\CB_x)_\phi$
is an $H_{4n}(Z)$-submodule of $H_{2d_x}(\CB_x)$.

The next theorem is proved directly in
\cite{kazhdanlusztig:topological} and \cite[\S3.5]
{chrissginzburg:representation}. It also follows from the
sheaf-theoretic approach to Borel-Moore homology described below.

\begin{theorem}\label{isotypic}
  There is an isomorphism of $H_{4n}(Z)$-bimodules,
  \[
  (H_{2d_x} (\CB_x) \otimes H_{2d_x} (\CB_x))^{C(x)} \cong
  \bigoplus_{\phi\in \widehat{C(x)}} \End_{\BBQ}( H_{2d_x}
  (\CB_x)_\phi).
  \]
  Moreover, $H_{2d_x}(\CB_x)_\phi$ is a simple $\gr\,H_{4n}(Z)$-module
  for every $\phi$ in $\widehat{C(x)}$ and the decomposition
  \[
  \gr\, H_{4n}(Z) \cong \bigoplus_{x\in \FS} \bigoplus_{\phi\in
    \widehat{C(x)}} \End_{\BBQ}( H_{2d_x}(\CB_x)_\phi)
  \]
  is a decomposition of $H_{4n}(Z)$ into minimal two-sided ideals.
\end{theorem}

Now that we have described the Wedderburn decomposition of $H_{4n}(Z)$
and given an explicitly construction of the irreducible
representations of $W$, we take up the question of finding formulas
for the action of a simple reflection.

For $x$ in $\FN$, formulas for the action of a simple reflection on
the basis of $H_{2d_x}(\CB_x)$ given by the irreducible components
were first given by Hotta and then refined by Borho, Brylinski, and
MacPherson (see \cite{hotta:local} and
\cite[\S4.14]{borhobrylinskimacpherson:nilpotent}). Analogous formulas
for the action of a simple reflection on $H_{4n}(Z)$ have been given
by Hinich and Joseph \cite[\S5]{hinichjoseph:orbital}. The first two
parts of the next theorem may be recovered from the more general (and
more complicated) argument in \cite[\S5]{douglassroehrle:homology}.

\begin{theorem}
  Suppose that $s$ is a simple reflection in $W$ and $w$ is in $W$. 
  \begin{itemize}
  \item[(a)] $\lambda_s= [\overline{Z_s}] +1$.
  \item[(b)] If $sw<w$, then $[\overline{Z_s}] *[\Zwbar]= -2
    [\Zwbar]$. 
  \item[(c)] If $sw>w$, then there is a subset $F_{s,w}$ of $\{\, x\in
    W\mid x<w,\, sx<x \,\}$ so that $[\overline{Z_s}] *[\Zwbar]=
    [\overline{Z_{sw}}] + \sum_{x\in F_{s,w}} n_x [\overline{Z_x}]$
    with $n_x>0$.
  \end{itemize}
\end{theorem}

Using this result, Hinich and Joseph \cite[Theorem
5.5]{hinichjoseph:orbital} prove a result analogous to Proposition
\ref{ideal} for right Steinberg cells. Recall that for $w$ in $W$ we
have defined $\FV_r(w)= \overline{B \fu_w} \cap \FO$ when $w$ is in
$W_\FO$. For an orbital variety $\FV$, define $W_{\FVbar} = \{\, y\in
W\mid \FV_r(y)\subseteq \overline{ \FV} \,\}$.

\begin{theorem}
  For $w$ in $W$, the smallest subset, $S$, of $W$ with the property
  that $[\Zwbar] *\lambda_y$ is in the span of $\{\, [\overline{Z_x}]
  \mid x\in S\,\}$ for all $y$ in $W$ is $\overline{\FV_r(w)}$. In
  particular, if $\FV$ is any orbital variety, then the span of $\{\,
  [\overline{Z_x}] \mid x\in W_{\FVbar}\,\}$ is a right ideal in
  $H_{4n}(Z)$.
\end{theorem}

\subsection{Sheaf-theoretic decomposition of $H_{4n}(Z)$ and
  $H_i(\CB_x)$} \label{s3.6}

For a variety $X$, the $\BBQ$-vector space $H_i(X)$ has more a
sophisticated alternate description in terms of sheaf cohomology (see
\cite[\S8.3] {chrissginzburg:representation}). The properties of
sheaves and perverse sheaves we use in this section may be found in
\cite[Chapter 2,3]{kashiwaraschapira:sheaves}, \cite{dimca:sheaves}
and \cite{borel:intersection}.

Let $D(X)$ denote the full subcategory of the derived category of
sheaves of $\BBQ$-vector spaces on $X$ consisting of complexes with
bounded, constructible cohomology. If $f\colon X\to Y$ is a morphism,
then there are functors 
\[
Rf_*\colon D(X)\to D(Y),\ Rf_!\colon D(X)\to D(Y), \ f^*\colon D(Y)\to
D(X), \ \text{and\ } f^!\colon D(Y)\to D(X).
\]
The pair of functors $(f^*,Rf_*)$ is an adjoint pair, as is $(Rf_!,
f^!)$. If $f$ is proper, then $Rf_!=Rf_*$ and if $f$ is smooth, then
$f^!=f^*[2\dim X]$. 

We consider the constant sheaf, $\BBQ_X$, as a complex in $D(X)$
concentrated in degree zero. The dualizing sheaf, $\BBD_X$, of $X$ is
defined by $\BBD_X=a_X^!\BBQ_{\{ \operatorname{pt} \}}$, where
$a_X\colon X\to \{ \operatorname{pt} \}$. If $X$ is a rational
homology manifold, in particular, if $X$ is smooth, then $\BBD_X \cong
\BBQ_X[2\dim X]$ in $D(X)$. It follows from the definitions and
because $f^*$ and $f^!$ are functors that if $f\colon X\to Y$, then
\begin{equation}
  \label{pullback}
  \BBQ_X\cong f^* \BBQ_Y \quad\text{and}\quad \BBD_X\cong f^! \BBD_Y  
\end{equation}
in $D(X)$.

The cohomology and Borel-Moore homology of $X$ have very convenient
descriptions in sheaf-theoretic terms:
\begin{equation}
  \label{Hdef}
  H^i(X) \cong \Ext^i_{D(X)} \left( \BBQ_X, \BBQ_X \right) \quad
  \text{and}\quad H_i(X) \cong \Ext^{-i}_{D(X)} \left( \BBQ_X, \BBD_X
  \right)   
\end{equation}
where for $\CF$ and $\CG$ in $D(X)$, $\Ext_{D(X)}^i (\CF, \CG)=
\Hom_{D(X)} (\CF, \CG[i])$.

Now suppose that $f_i\colon M_i\to N$ is a proper morphism for
$i=1,2,3$ and that $d_2=\dim M_2$. In contrast to our assumptions in
the convolution setup from \S\ref{s3.1} where the $M_i$ were assumed
to be smooth, in the following computation we assume only that $M_2$
is a rational homology manifold.  Consider the cartesian diagram
\[
\xymatrix{
  M_1\times_N M_2 \ar[r]^{f_{1,2}} \ar[d]_{\delta_1} &N
  \ar[d]^{\delta} \\ 
  M_1\times M_2 \ar[r]_{f_1\times f_2} & N\times N}
\]
where $f_{1,2}$ is the induced map. Using the argument in
\cite[\S8.6]{chrissginzburg:representation}, we have isomorphisms
\begin{align*}
  H_i(M_1\times_N M_2) &\cong\Ext^{-i}_{D(M_1\times_N M_2)} (\BBQ_{
    M_1 \times_N M_2}, \BBD_{M_1 \times_N M_2} )&
  \text{(\ref{Hdef})}\\
  &\cong \Ext^{-i}_{D(M_1\times_N M_2)} (f_{1,2}^* \BBQ_{ N},
  \delta_1^! \BBD_{M_1 \times M_2} )& \text{(\ref{pullback})}\\
  &\cong \Ext^{-i}_{D(N)} ( \BBQ_{ N}, R(f_{1,2})_*\delta_1^!
  \BBD_{M_1 \times M_2} )& \text{(adjunction)}\\
  &\cong \Ext^{-i}_{D(N)} ( \BBQ_{ N}, \delta^! R(f_1\times f_2)_*
  \BBD_{M_1 \times M_2} )& \text{(base change)}\\
  &\cong \Ext^{-i}_{D(N)} ( \BBQ_{ N}, \delta^! ( R(f_1)_*
  \BBD_{M_1} \boxtimes R(f_2)_*  \BBD_{M_2} ))& \text{(K\"unneth)}\\
  &\cong \Ext^{-i}_{D(N)} ( \BBQ_{ N}, \mathcal Hom( R(f_1)_*
  \BBQ_{M_1}, R(f_2)_* \BBD_{M_2} ))& \text{(\cite[10.25]
    {borel:intersection})}\\
  &\cong \Ext^{-i}_{D(N)} ( \BBQ_{ N}, \mathcal Hom( R(f_1)_*
  \BBQ_{M_1}, R(f_2)_* \BBQ_{M_2}[2d_2] ))& (\BBD_{M_2} \cong
  \BBQ_{M_2}[2d_2]) \\ 
  &\cong \Ext^{2d_2-i}_{D(N)} ( \BBQ_{ N}, \mathcal Hom( R(f_1)_*
  \BBQ_{M_1}, R(f_2)_* \BBQ_{M_2} ))&  \\
  &\cong \Ext^{2d_2-i}_{D(N)} ( R(f_1)_* \BBQ_{M_1}, R(f_2)_*
  \BBQ_{M_2} ). &
\end{align*}
Let $\epsilon_{1,2}$ denote the composition of the above isomorphisms,
so
\begin{equation}
  \label{fibre}
  \epsilon_{1,2} \colon H_i(M_1\times_N M_2)\xrightarrow{\cong}
  \Ext^{2d_2-i}_{D(N)} ( R(f_1)_* \BBQ_{M_1}, R(f_2)_* \BBQ_{M_2} ).  
\end{equation}
Chriss and Ginzburg \cite[\S8.6]{chrissginzburg:representation} have
shown that the isomorphisms $\epsilon_{1,2}$ intertwine the
convolution product on the left with the Yoneda product (composition
of morphisms) on the right: given $c$ in $H_i(M_1\times_N M_2)$
and $d$ in $H_j(M_2\times_N M_3)$, we have $\epsilon_{1,3}(c*d)=
\epsilon_{2,3}(d) \circ \epsilon_{1,2}(c)$. 

We may apply the computation in equation (\ref{fibre}) to $H_*(Z)$. We
have seen that $Z\cong \FNt \times_\FN \FNt$ and so
\[
H_i(Z) \cong \Ext^{4n-i} _{D(\FN)} ( R\mu_* \BBQ_{\FNt}, R\mu_*
\BBQ_{\FNt}).
\]
In particular, taking $i=4n$, we conclude that are algebra isomorphisms
\[
\BBQ[W] \cong H_{4n}(Z) \cong \End _{D(\FN)} ( R\mu_*
\BBQ_{\FNt})^{\op}.
\]

The category $D(\FN)$ is a triangulated category. It contains a full,
abelian subcategory, denoted by $\CM(\FN)$, consisting of ``perverse
sheaves on $\FN$'' (with respect to the middle perversity). It follows
from the Decomposition Theorem of Beilinson, Bernstein, and Deligne
\cite[\S5]{beilinsonbernsteindeligne:faisceaux} that the complex
$R\mu_* \BBQ_\FNt$ is a semisimple object in $\CM(\FN)$.

The simple objects in $\CM(\FN)$ have a geometric description. Suppose
$X$ is a smooth, locally closed subvariety of $\FN$ with codimension
$d$, $i\colon X\to \FN$ is the inclusion, and $L$ is an irreducible
local system on $X$. Let $\IC(X, L)$ denote the intersection complex
of Goresky and MacPherson \cite[\S3]
{goreskymacpherson:intersectionII}.  Then $i_* \IC(\overline{X} ,
L)[-2d]$ is a simple object in $\CM(\FN)$ and every simple object
arises in this way. In addition to the original sources,
\cite{beilinsonbernsteindeligne:faisceaux} and
\cite{goreskymacpherson:intersectionII}, we refer the reader to
\cite[\S3]{shoji:geometry} and
\cite[\S8.4]{chrissginzburg:representation} for short introductions to
the theory of intersection complexes and perverse sheaves and to
\cite{borel:intersection} and \cite{dimca:sheaves} for more thorough
expositions.

Returning to $R\mu_* \BBQ_\FNt$, Borho and MacPherson
\cite{borhomacpherson:weyl} have shown that its decomposition into
simple perverse sheaves is given by
\begin{equation}
  \label{eq:pv}
  R\mu_* \BBQ_\FNt \cong \bigoplus_{x, \phi} j_*^x \IC( \overline{Gx},
  L_\phi)[-2d_x] ^{n_{x,\phi}} 
\end{equation}
where $x$ runs over the set of orbit representatives $\FS$ in $\FN$,
and for each $x$, $j^x\colon \overline{Gx} \to \FN$ is the inclusion,
$\phi$ is in $\widehat{C(x)}$, $L_\phi$ is the local system on $Gx$
corresponding to $\phi$, and $n_{x,\phi}$ is a non-negative integer.

Define $\IC_{x,\phi}= j_*^x \IC( \overline{Gx}, L_\phi)$. Then
$\IC_{x,\phi}[-2d_x]$ is a simple object in $\CM(\FN)$. It follows
from the computation of the groups $C(x)$ that $\End_{D(\FN)}
(\IC_{x,\phi}) \cong \BBQ$. Therefore,
\begin{equation}
  \label{h4nz}
  \begin{aligned}
    H_{4n}(Z) &\cong \End _{D(\FN)} ( R\mu_* \BBQ_{\FNt})^{\op} \\
    &\cong \End _{D(\FN)} ( \oplus_{x, \phi} \IC_{x,\phi}[-2d_x]
    ^{n_{x,\phi}}  )^{\op} \\
    & \cong \bigoplus_{x, \phi} \End _{D(\FN)} ( \IC_{x,\phi}
    ^{n_{x,\phi}} )^{\op}\\
    & \cong \bigoplus_{x, \phi} M_{n_{x, \phi}} \left( \End _{D(\FN)}
      (
      \IC_{x,\phi} ) \right)^{\op}\\
    & \cong \bigoplus_{x, \phi} M_{n_{x, \phi}} \left(\BBQ
    \right)^{\op}.
  \end{aligned}
\end{equation}
This is a decomposition of $H_{4n}(Z)$ as a direct sum of matrix rings
and hence is the Wedderburn decomposition of $H_{4n}(Z)$.

Suppose now that $\FO$ is a $G$-orbit in $\FN$ and $x$ is in $\FO$.
It is straightforward to check that
\[
H_\FO \cong \bigoplus_{\phi\in \widehat{ C(x)}} \End_{D(\FN)} (
(\IC_{x,\phi}) ^{n_{x,\phi}}) \cong \bigoplus_{\phi\in \widehat{ C(x)}}
M_{n_{x, \phi}}(\BBQ).
\]
As in Proposition \ref{ideal}, this is the decomposition of $H_\FO$
into minimal two-sided ideals.

For a second application of (\ref{fibre}), let $i_x\colon \{x\} \to
\FN$ denote the inclusion. Then $\CB_x \cong \FNt \times_{\FN} \{x\} $
and so
\begin{align*}
  H_i(\CB_x) &\cong \Ext_{D(\FN)} ^{-i}(R\mu_* \BBQ_{\FNt}, R(i_x)_*
  \BBQ_{\{x\}} )\\
  &\cong \bigoplus_{y, \psi} \Ext_{D(\FN)} ^{-i}( \IC_{y, \psi}[-2d_y]
  ^{n_{y,\psi}}, R(i_x)_* \BBQ_{\{x\}} )\\
  &\cong \bigoplus_{y, \psi} \Ext_{D(\FN)} ^{2d_y-i}(
  \IC_{y, \psi} ^{n_{y,\psi}}, R(i_x)_* \BBQ_{\{x\}} )\\
  &\cong \bigoplus_{y, \psi} \left( V_{y, \psi} \otimes \Ext_{D(\FN)}
    ^{2d_y-i}( \IC_{y, \psi}, R(i_x)_* \BBQ_{\{x\}}) \right)
\end{align*}
where $V_{y,\psi}$ is an $n_{y, \psi}$-dimensional vector space.
Because $\BBQ[W] \cong H_{4n}(Z) \cong \End _{D(\FN)} ( R\mu_*
\BBQ_{\FNt})$ acts by permuting the simple summands, it follows from
(\ref{h4nz}) that each $V_{y, \psi}$ affords an irreducible
representation of $W$ and that $\Ext_{D(\FN)} ^{2d_y-i}( \IC_{y, \psi},
R(i_x)_* \BBQ_{\{x\}})$ records the multiplicity of $V_{y, \psi}$ in
$H_i(\CB_x)$. Using that $i_x^*$ is left adjoint to $R(i_x)_*$,
denoting the stalk of $\IC_{y, \psi}$ at $x$ by $(\IC_{y, \psi})_x$, and
setting $m_{y,\psi}^{x,i} = \dim \Ext_{D(\FN)} ^{2d_y-i}( \IC_{y,
  \psi}, R(i_x)_* \BBQ_{\{x\}})$, we obtain the decomposition of
$H_i(\CB_x)$ into irreducible representations of $W$:
\[
H_i(\CB_x) \cong \bigoplus_{y, \psi} \left( V_{y, \psi} \otimes
  \Ext_{D(\{x\} )} ^{2d_y-i}( (\IC_{y, \psi})_x, \BBQ_{\{x\}}) \right)
\cong \bigoplus_{y, \psi} V_{y, \psi} ^{m_{y,\psi}^{x,i}}.
\]

In the next subsection we apply (\ref{fibre}) to compute the
Borel-Moore homology of some generalized Steinberg varieties.

\subsection{Borel-Moore homology of generalized Steinberg varieties}
\label{s3.7} 

Recall from \S\ref{s2.4} the generalized Steinberg variety 
\[
\XPQ= \{\, (x, P',Q')\in \FN \times \CP \times \CQ\mid x\in \fp' \cap
\fq'\,\} \cong \FNt^{\CP} \times_{\FN}\FNt^{\CQ}
\]
where $\FNt^{\CP}= \{\, (x, P')\in \FN\times \CP \mid x\in \fp'\,\}$,
$\xi^{\CP}\colon \FNt^{\CP} \to \FN$ is projection on the first
factor, and $\FNt^{\CQ}$ and $\xi^{\CQ}$ are defined similarly.
Recall also that $\eta\colon Z\to \XPQ$ is a proper, $G$-equivariant
surjection. The main result of \cite[Theorem 4.4]
{douglassroehrle:homology}, which is proved using the constructions in
the last subsection, is the following theorem describing the
Borel-Moore homology of $\XPQ$.

\begin{theorem}\label{ave}
  Consider $H_{4n}(Z)$ as a $W\times W$-module using the isomorphism
  $H_{4n}(Z) \cong \BBQ[W]$. Then there is an isomorphism
  $\alpha\colon H_*(\XPQ) \xrightarrow{\cong} H_*(Z)^{W_P\times W_Q}$
  so that the composition $\alpha \circ \eta_*\colon H_*(Z) \to
  H_*(Z)^{W_P\times W_Q}$ is the averaging map.
\end{theorem}

As a special case of the theorem, if we let $e_P$ (resp.~$e_Q$) denote
the primitive idempotent in $\BBQ[W_P]$ (resp.~$\BBQ[W_Q]$)
corresponding to the trivial representation, then   
\begin{equation}
  \label{X}
  H_{4n}(\XPQ) \cong e_P \BBQ[W] e_Q.  
\end{equation}

Next recall the generalized Steinberg variety $\XPQoo \cong T^*\CP
\times_\FN T^* \CQ$. Set $m= \dim P/B +\dim Q/B$.  Let $\epsilon_P$
(resp.~$\epsilon_Q$) denote the primitive idempotent in $\BBQ[W_P]$
(resp.~$\BBQ[W_Q]$) corresponding to the sign representation. Then
$\dim \XPQoo= 4n-2m$ and it is shown in
\cite[\S5]{douglassroehrle:homology} that 
\begin{equation}
  \label{X0}
  H_{4n-2m}(\XPQoo) \cong \epsilon_P \BBQ[W] \epsilon_Q.  
\end{equation}

Now suppose that $c$ is a Levi class function on $\CP$. Let $L$ be a
Levi subgroup of $P$ and choose $x$ in $c(P)\cap \fl$. Then we may
consider the Springer representation of $W_P$ on
$H_{2d^L_x}(\CB^L_x)^{C_L(x)}$ where $C_L(x)$ is the component group
of $Z_L(x)$, $\CB^L_x$ is the variety of Borel subalgebras of $\fl$
that contain $x$, and $d^L_x= \dim \CB^L_x$. This is an irreducible
representation of $W_P$. Let $f_P$ denote a primitive idempotent in
$\BBQ[W_P]$ so that $\BBQ[W_P] f_P\cong H_{2d^L_x}(\CB^L_x)^{C_L(x)}$.
Set $\delta^{\CP, \CQ}_{c,d}= \frac 12\left(\dim c(P)+ \dim \fu_P+
  \dim d(Q)+ \dim \fu_Q\right)$.  Then it is shown in \cite[Corollary
2.6]{douglassroehrle:geometry} that $\dim \XPQcd\leq \delta^{\CP,
  \CQ}_{c,d}$. Generalizing the computations (\ref{X}) and (\ref{X0}),
we conjecture that the following statement is true.
\begin{conjecture}\label{conj}
  With the notation above, $H_{\delta^{\CP, \CQ}_{c,d}}(\XPQcd) \cong
  f_P \BBQ[W] f_Q$.
\end{conjecture}

The Borel-Moore homology of $\XPQ$ may also be computed using the
sheaf theoretic methods in the last subsection. We have $\XPQ \cong
\FNt^{\CP} \times_{\FN}\FNt^{\CQ}$ and Borho and MacPherson
\cite[2.11]{borhomacpherson:partial} have shown that $\FNt^{\CP}$ and
$\FNt^{\CQ}$ are rational homology manifolds. Therefore, as in
(\ref{fibre}):
\[
H_i(\XPQ) \cong \Ext_{D(\FN)} ^{4n-i} ( R\xi^{\CP}_*
\BBQ_{ \FNt^{\CP}}, R\xi^{\CQ}_* \BBQ_{ \FNt^{\CQ}}).
\]

Borho and MacPherson \cite[2.11]{borhomacpherson:partial} have also
shown that $R\xi^{\CP}_* \BBQ_{ \FNt^{\CP}}$ is a semisimple object in
$\CM(\FN)$ and described its decomposition into simple perverse
sheaves:
\[
R\xi^{\CP}_* \BBQ_{ \FNt^{\CP}} \cong \bigoplus_{(x, \phi)}
\IC_{x,\phi} [-2d_x]^{n^\CP_{x, \phi}},
\]
where the sum is over pairs $(x,\phi)$ as in equation (\ref{eq:pv}),
and $n^\CP_{x, \phi}$ is the multiplicity of the irreducible
representation $H_{2d_x}(\CB_x)_\phi$ of $W$ in the induced
representation $\Ind_{W_P}^W (1_{W_P})$. Thus,
\[
H_i(\XPQ) \cong \bigoplus_{x, \phi} \bigoplus_{y, \psi}
\Ext_{D(\FN)}^{4n-i} \left( \IC_{x,\phi}[-2d_x] ^{n^\CP_{x, \phi}} ,
  \IC_{y,\psi} [-2d_y]^{n^\CQ_{y, \psi}} \right)
\]
and so
\begin{equation} \label{4nxpq} 
  H_{4n}(\XPQ) \cong \bigoplus_{x, \phi}
  \bigoplus_{y, \psi} \Hom_{D(\FN)} \left( \IC_{x,\phi}
    [-2d_x]^{n^\CP_{x, \phi}} , \IC_{y,\psi} [-2d_y] ^{n^\CQ_{y,
        \psi}} \right) \cong \bigoplus_{x, \phi} M_{n^\CQ_{x, \phi},
    n^\CP_{x, \phi}} (\BBQ).
\end{equation}

Using the fact that $n^\CP_{x, \phi}$ is the multiplicity of the
irreducible representation $H_{2d_x}(\CB_x)_\phi$ of $W$ in the
induced representation $\Ind_{W_P}^W (1_{W_P})$, we see that
(\ref{4nxpq}) is consistent with (\ref{X}).

\section{Equivariant $K$-theory}

Certainly the most important result to date involving the Steinberg
variety is its application by Kazhdan and Lusztig to the Langlands
program \cite{kazhdanlusztig:langlands}. They show that the
equivariant $K$-theory of $Z$ is isomorphic to the two-sided regular
representation of the extended, affine Hecke algebra $\CH$. They then
use this representation of $\CH$ to classify simple $\CH$-modules and
hence to classify representations of ${}^LG( \BBQ_p)$ containing a
vector fixed by an Iwahori subgroup, where ${}^LG( \BBQ_p)$ is the
group of $\BBQ_p$-points of the Langlands dual of $G$. As with homology,
Chriss and Ginzburg have applied the convolution product formalism to
the equivariant $K$-theory of $Z$ and recast Kazhdan and Lusztig's
results as an algebra isomorphism.

Recall we are assuming that $G$ is simply connected. In this section
we describe the isomorphism $\CH\cong K^{\Gbar}(Z)$, where $\Gbar=
G\times \BBC^*$, and we give some applications to the study of
nilpotent orbits. We emphasize in particular the relationship between
nilpotent orbits, Kazhdan-Lusztig theory for the extended, affine Weyl
group, and (generalized) Steinberg varieties.

\subsection{The generic, extended, affine Hecke algebra} \label{s4.1}
We begin by describing the Bern\-stein-Zelevinski presentation of the
extended, affine Hecke algebra following the construction in
\cite{lusztig:affine}.

Let $v$ be an indeterminate and set $A=\BBZ[v, v\inverse]$. The ring
$A$ is the base ring of scalars for most of the constructions in this
section. 

Let $X(T)$ denote the character group of $T$. Since $G$ is simply
connected, $X(T)$ is the weight lattice of $G$. Define $X^+$ to be the
set of dominant weights with respect to the base of the root system of
$(G,T)$ determined by $B$. The \emph{extended, affine Weyl group} is
$W_e= X(T)\rtimes W$. 

There is a ``length function'' $\ell$ on $W_e$ that extends the usual
length function on $W$. The \emph{braid group of $W_e$} is the group
$\mathcal Br$, with generators $\{\, T_x\mid x\in W_e\,\}$ and
relations $T_xT_{x'}= T_{xx'}$ if $\ell(x)+ \ell(x')= \ell(xx')$. The
\emph{generic, extended, affine Hecke algebra,} $\CH$, is the quotient
of the group algebra $A[\mathcal Br]$ by the two-sided ideal generated
by the elements $(T_s+1) (T_s-v^2)$, where $s$ runs through the simple
reflections in $W$.

Let ${}^LG$ denote the Langlands dual of $G$, so ${}^LG$ is an adjoint
group. Let ${}^LG_p$ denote the algebraic group over $\BBQ_p$ with the
same type as ${}^LG$. Suppose that $I$ is an Iwahori subgroup of
${}^LG_p$ and let $\BBC[I\backslash {}^LG_p/I]$ denote the space of
all compactly supported functions ${}^LG_p\to \BBC$ that are constant
on $(I,I)$-double cosets. Consider $\BBC$ as an $A$-module via the
specialization $A\to \BBC$ with $v\mapsto \sqrt p$. The following
theorem, due to Iwahori and Matsumoto \cite[\S3]
{iwahorimatsumoto:bruhat}, relates $\CH$ to representations of
${}^LG_p$ containing an $I$-fixed vector.

\begin{theorem}
  The $(I,I)$-double cosets of ${}^LG_p$ are parametrized by
  $W_e$. Moreover, if $I_x$ is the double coset indexed by $x$ in
  $W_e$, then the map which sends $T_x$ to the characteristic function
  of $I_x$ extends to an algebra isomorphism
  \[
  \BBC\otimes_A \CH \cong \BBC[I\backslash {}^LG_p/I].
  \]
\end{theorem}

The algebra $\CH$ has a factorization (as a tensor product) analogous
to the factorization $W_e =X(T)\rtimes W$. Given $\lambda$ in $X(T)$
one can write $\lambda= \lambda_1- \lambda_2$ where $\lambda_1$ and
$\lambda_2$ are in $X^+$.  Define $E^\lambda$ in $\CH$ to be the image
of $v^{\ell( \lambda_1 -\lambda_2)} T_\lambda$. For $x$ in $W_e$,
denote the image of $T_x$ in $\CH$ again by $T_x$. Let $\CH_W$ denote
the Iwahori-Hecke algebra of $W$ (an $A$-algebra) with standard basis
$\{ t_w\mid w\in W\,\}$.  Lusztig \cite[\S2] {lusztig:affine} has
proved the following theorem.

\begin{theorem}
  With the notation above we have:
  \begin{itemize}
  \item[(a)] $E^\lambda$ does not depend on the choice of $\lambda_1$
    and $\lambda_2$.
  \item[(b)] The mapping $A[X(T)] \otimes _A \CH_W \to \CH$ defined by
    $\lambda \otimes t_w\mapsto E^\lambda T_w$ is an isomorphism of
    $A$-modules.
  \item[(c)] For $\lambda$ and $\lambda'$ in $X$ we have $E^\lambda
    E^{\lambda'}= E^{\lambda +\lambda'}$ and so the subspace of $\CH$
    spanned by $\{\, E^\lambda \mid \lambda\in X\,\}$ is a subalgebra
    isomorphic to $A[X(T)]$.
  \item[(d)] The center of $\CH$ is isomorphic to $A[X(T)]^W$ via the
    isomorphism in (c).
  \item[(e)] The subspace of $\CH$ spanned by $\{\, T_w\mid w\in
    W\,\}$ is a subalgebra isomorphic to $\CH_W$.
  \end{itemize}
\end{theorem}

Using parts (b) and (d) of the theorem, we identify $A[X(T)]$ with the
subalgebra of $\CH$ spanned by $\{\, E^\lambda \mid \lambda\in X\,\}$,
and $A[X(T)]^W$ with the center of $\CH$.

\subsection{Equivariant $K$-theory and convolution}
Two introductory references for the notions from equivariant
$K$-theory we use are \cite[Chapter 2]
{borhobrylinskimacpherson:nilpotent} and \cite[Chapter 5]
{chrissginzburg:representation}.

For a variety $X$, let $\mathcal Coh(X)$ denote the category of
coherent $\CO_X$-modules. Suppose that $H$ is a linear algebraic group
acting on $X$. Let $a\colon H\times X\to X$ be the action morphism and
$p\colon H\times X\to X$ be the projection. An \emph{$H$-equivariant
  coherent $\CO_X$-module} is a pair $(\CM, i)$, where $\CM$ is a
coherent $\CO_X$-module and $i\colon a^*\CM\xrightarrow{\sim} p^* \CM$
is an isomorphism satisfying several conditions (see \cite[\S5.1]
{chrissginzburg:representation} for the precise definition). With the
obvious notion of morphism, $H$-equivariant $\CO_X$-modules form an
abelian category denoted by $\mathcal Coh^H(X)$. The Grothendieck
group of $\mathcal Coh^H(X)$ is denoted by $K^H(X)$ and is called the
\emph{$H$-equivariant $K$-group of $X$.}

If $X=\{ \operatorname{pt}\}$ is a point, then $K^H (\operatorname
{pt})\cong R(H)$ is the representation ring of $H$.  For any $X$,
$K^H(X)$ is naturally an $R(H)$-module. If $H$ is the trivial group,
then $\mathcal Coh^H(X)= \mathcal Coh(X)$ and $K^H(X)=K(X)$ is the
Grothendieck group of the category of coherent $\CO_X$-modules.

As with Borel-Moore homology, equivariant $K$-theory is a bivariant
theory in the sense of Fulton and MacPherson
\cite{fultonmacpherson:categorical}: Suppose that $X$ and $Y$ are
$H$-varieties and that $f\colon X\to Y$ is an $H$-equivariant
morphism. If $f$ is proper, there is a direct image map in equivariant
$K$-theory, $f_*\colon K^H(X) \to K^H(Y)$, and if $f$ is smooth there
is a pullback map $f^*\colon K^H(Y) \to K^H(X)$ in equivariant
$K$-theory. Moreover, if $X$ is smooth and $A$ and $B$ are closed,
$H$-stable subvarieties of $X$, there is an intersection pairing
$\cap\colon K^H(A) \times K^H(B)\to K^H(A\cap B)$ (called a
Tor-product in \cite[\S6.4] {lusztig:bases}).  This pairing depends on
$(X,A,B)$. Thus, we may apply the convolution product construction
from \S\ref{s3.1} in the equivariant $K$-theory setting.

In more detail, suppose that for $i=1,2,3$, $M_i$ is a smooth variety
with an algebraic action of $H$ and $f_i\colon M_i\to N$ is a proper,
$H$-equivariant morphism. Suppose that for $1\leq i<j\leq 3$,
$Z_{i,j}$ is a closed, $H$-stable subvariety of $M_i\times M_j$ and
that $p_{1,3}\colon p_{1,2}\inverse(Z_{1,2}) \cap p_{2,3} \inverse
(Z_{2,3}) \to Z_{1,3}$ is a proper morphism. Then as in \S\ref{s3.1},
the formula $c*d= (p_{1,3})_*\left( p_{1,2}^*(c) \cap p_{2,3}^*(d)
\right)$, where $\cap$ is the intersection pairing determined by the
subsets $Z_{1,2} \times M_3$ and $M_1 \times Z_{2,3}$ of $M_1\times
M_2\times M_3$, defines an associative convolution product,
$K^H(Z_{1,2}) \otimes K^H(Z_{2,3})\xrightarrow{*} K^H (Z_{1,3})$.

In particular, the convolution product defines a ring structure on
$K^G(Z)$. It is shown in \cite[Theorem
7.2.2]{chrissginzburg:representation} that with this ring structure,
$K^G(Z)$ is isomorphic to the group ring $\BBZ[W_e]$. In the next
subsection we describe a more general result with $\BBZ[W_e]$ replaced
by $\CH$ and $G$ replaced by $G\times \BBC^*$, where $\BBC^*$ denote
the multiplicative group of non-zero complex numbers. 

The variable, $v$, in the definition of $\CH$ is given a geometric
meaning using the isomorphism $X(\BBC^*) \cong \BBZ$. Let $1_{\BBC^*}$
denote the trivial representation of $\BBC^*$. Then the rule $v\mapsto
1_{\BBC^*}$ extends to a ring isomorphism $\BBZ[v,v\inverse] \cong
R(\BBC^*)$. For the rest of this paper we will use this isomorphism to
identify $A=\BBZ[v,v\inverse]$ and $R(\BBC^*)$.

\subsection{The Kazhdan-Lusztig isomorphism} \label{s4.3}
To streamline the notation, set $\Gbar= G\times \BBC^*$. Then
$R(\Gbar) \cong R(G) \otimes_{\BBZ} R(\BBC^*) \cong R(G)
\otimes_{\BBZ} A= R(G)[v,v\inverse]$. 

Similarly, for a closed subgroup, $H$, of $G$, we denote the subgroup
$H\times \BBC^*$ of $\Gbar$ by $\overline{H}$. In particular, $\Tbar=
T\times \BBC^*$ and $\overline{B}= B\times \BBC^*$. In the remainder
of this paper we will never need to consider the closure of a subgroup
of $G$ and so this notation should not lead to any confusion.

Define a $\BBC^*$-action on $\fg$ by $(\xi, x)\mapsto \xi^{2} x$. We
consider $\CB$ as a $\BBC^*$-set with the trivial action. Then the
action of $G$ on $\FNt$ and $Z$ extends to an action of $\Gbar$ on
$\FNt$ and $Z$, and $\mu_z$ and $\mu$ are $\Gbar$-equivariant.

Recall from \S\ref{s4.1} that we are viewing the group algebra
$A[X(T)]$ as a subspace of $\CH$, and that the center of $\CH$ is
$Z(\CH)=A[X(T)]^W$. Using the identification $A=R(\BBC^*)$, we may
begin to interpret subspaces of $\CH$ in $K$-theoretic terms:
\[
K^{\Gbar}(\{ \operatorname{pt} \}) \cong R(\Gbar) \cong R(G)\otimes
R(\BBC^*) \cong R(G)[v,v\inverse] \cong A [X(T)]^W= Z(\CH).
\]

Recall that the ``diagonal'' subvariety, $Z_1$, of the Steinberg
variety is defined by $Z_1=\{\, (x, B', B')\in \FN\times \CB\times \CB
\mid x\in \fb'\,\}$.  For suitable choices of $f_i\colon M_i\to N$ and
$Z_{i,j}$, and using the embedding $A\subseteq R(\Gbar)$, the
convolution product induces various $A$-linear maps:
\begin{itemize}
\item[(1)] $K^{\Gbar}(Z) \times K^{\Gbar}(Z) \xrightarrow{*}
  K^{\Gbar}(Z)$; with this multiplication, $K^{\Gbar}(Z)$ is an
  $A$-algebra.
\item[(2)] $K^{\Gbar}(Z_1) \times K^{\Gbar}(Z_1) \xrightarrow{*}
  K^{\Gbar}(Z_1)$; with this multiplication, $K^{\Gbar}(Z_1)$ is a
  commutative $A$-algebra.
\item[(3)] $K^{\Gbar}(Z) \times K^{\Gbar}(\FNt \times \CB)
  \xrightarrow{*} K^{\Gbar}( \FNt \times \CB)$; this defines a left
  $K^{\Gbar}(Z)$-module structure on $K^{\Gbar}(\FNt \times \CB)$.
\end{itemize}

The group $K^{\Gbar}(Z_1)$ has a well-known description. First, the
rule $(x, B')\mapsto (x, B', B')$ defines a $\Gbar$-equivariant
isomorphism between $\FNt$ and $Z_1$ and hence an isomorphism
$K^{\Gbar}( Z_1) \cong K^{\Gbar}( \FNt)$. Second, the projection
$\FNt\to \CB$ is a vector bundle and so, using the Thom isomorphism in
equivariant $K$-theory \cite[\S5.4]{chrissginzburg:representation}, we
have $K^{\Gbar} (\FNt) \cong K^{\Gbar} (\CB)$. Third, $\CB$ is
isomorphic to $G \times^{B} \{ \operatorname{pt} \}$ by a
$\Gbar$-equivariant isomorphism and so $K^{\Gbar}( \CB) \cong K^{
  \overline B} ( \{ \operatorname{pt} \}) \cong R(\overline B)$ by a
version of Frobenius reciprocity in equivariant $K$-theory
\cite[\S5.2.16]{chrissginzburg:representation}. Finally, since $U$ is
the unipotent radical of $B$, we have 
\[
R( \overline B)\cong R(\overline B/\overline U) \cong R(\Tbar) \cong
R(T)[v,v\inverse] \cong A[X(T)].
\]
Composing these isomorphisms, we get an isomorphism $K^{\Gbar}(Z_1)
\xrightarrow {\cong} A[X(T)]$, which is in fact an isomorphism of
$A$-algebras.

The inverse isomorphism $A[X(T)] \xrightarrow{\cong} K^{\Gbar}(Z_1)$
may be computed explicitly. Suppose that $\lambda$ is in $X(T)$. Then
$\lambda$ lifts to a representation of $\overline{B}$. Denote the
representation space by $\BBC_\lambda$. Then the projection morphism
$\Gbar\times^{\overline B} \BBC_\lambda \to \CB$ is a
$\Gbar$-equivariant line bundle on $\CB$. The sheaf of sections of
this line bundle is a $\Gbar$-equivariant, coherent sheaf of
$\CO_\CB$-modules that we will denote by $L_\lambda$. Pulling
$L_\lambda$ back first through the vector bundle projection $\FNt\to
\CB$ and then through the isomorphism $Z_1\cong \FNt$, we get a
$\Gbar$-equivariant, coherent sheaf of $\CO_{Z_1}$-modules we denote
by $\CL_\lambda$.

Let $i_1\colon Z_1\to Z$ be the inclusion. Define $e^\lambda =
(i_1)_*([\CL_\lambda])$ in $K^{\Gbar}(Z)$. Then $\lambda\mapsto
e^\lambda$ defines an $A$-linear map from $A[X(T)]$ to $K^{\Gbar}(Z)$.

A concentration theorem due to Thomason and the Cellular Fibration
Lemma of Chriss and Ginzburg can be used to prove the following
proposition (see \cite[6.2.7]{chrissginzburg:representation} and
\cite[7.15]{lusztig:bases}).

\begin{proposition}\label{faithful}
  The closed embeddings $i_1\colon Z_1\to Z$ and $j\colon Z\to
  \FNt\times \CB$ induce injective maps in equivariant $K$-theory,
  \[
  K^{\Gbar}(Z_1) \xrightarrow{(i_1)_*} K^{\Gbar}(Z) \xrightarrow{j_*}
  K^{\Gbar}( \FNt\times \CB).
  \]
  The map $(i_1)_*$ is an $A$-algebra monomorphism and the map $j_*$
  is a $K^{\Gbar}(Z)$-module monomorphism. In particular, $K^{\Gbar}(
  \FNt \times \CB)$ is a faithful $K^{\Gbar} (Z)$-module.
\end{proposition}

From the proposition and the isomorphism $K^{\Gbar} (\{ \operatorname
{pt} \}) \cong Z(\CH)$, we see that there is a commutative diagram of
$A$-algebras and $A$-algebra homomorphisms:
\[
\xymatrix{ Z(\CH) \ar@{^{(}->}[r] \ar[d]_{\cong} & A[X(T)]
  \ar@{^{(}->}[r] \ar[d]_{\cong}& \CH \\ K^{\Gbar} (\{
  \operatorname{pt} \}) \ar@{^{(}->}[r] & K^{\Gbar}(Z_1)
  \ar@{^{(}->}[r] &K^{\Gbar}(Z).}
\]
We will complete this diagram with an isomorphism of $A$-algebras
$K^{\Gbar}(Z) \cong \CH$ following the argument in
\cite[\S7]{lusztig:bases}.

Fix a simple reflection, $s$, in $W$. Then there is a simple root,
$\alpha$, in $X(T)$ and a corresponding cocharacter, $\alphacheck
\colon \BBC^*\to T$, so that if $\langle \,\cdot\, ,\,\cdot \,
\rangle$ is the pairing between characters and cocharacters of $T$,
then $\langle \alpha, \alphacheck \rangle =2$ and $s(\lambda)= \lambda
- \langle \lambda, \alphacheck \rangle \alpha$ for $\lambda$ in
$X(T)$. Choose a weight $\lambda'$ in $X(T)$ with $\langle \lambda',
\alphacheck \rangle = -1$ and set $\lambda''= - \lambda'-\alpha$. Then
$L_{\lambda'} \boxtimes L_{\lambda''}$ is in $\mathcal Coh^{\Gbar}(
\CB \times \CB)$. Lusztig \cite[7.19]{lusztig:bases} has shown that
the restriction of $L_{\lambda'} \boxtimes L_{\lambda''}$ to the
closed subvariety $\overline{ G(B, sBs)}$ does not depend on the
choice of $\lambda'$. Denote the restriction of $L_{\lambda'}
\boxtimes L_{\lambda''}$ to $\overline{ G(B, sBs)}$ by $\CL_s$.

It is easy to check that $Z_1\cap \overline{Z_s}= \{\,(x, gBg\inverse,
gBg\inverse)\in Z_1 \mid g\inverse x \in \fu_s\,\}$. It follows that
$\overline{Z_s}$ is smooth and that $\pi \colon \overline{Z_s} \to
\overline{ G(B, sBs)}$ is a vector bundle projection with fibre
$\fu_s$. Thus, there is a pullback map in equivariant $K$-theory,
$\pi^*\colon K^{\Gbar} \left( \overline{ G(B, sBs)} \right) \to
K^{\Gbar} \left (\overline{Z_s} \right)$, and so we may consider
$\pi^*( [\CL_s])$ in $K^{\Gbar} \left( \overline{Z_s} \right)$. Let
$i_s\colon \overline{Z_s} \to Z$ denote the inclusion. Then $i_s$ is a
closed embedding and so there is a direct image map $(i_s)_* \colon
K^{\Gbar} \left (\overline{Z_s} \right) \to K^{\Gbar}(Z)$. Define
$l_s= (i_s)_* \pi^*([\CL_s])$. Then $l_s$ is in $K^{\Gbar}(Z)$.

Lusztig \cite[7.24]{lusztig:bases} has proved the following lemma.

\begin{lemma}\label{contain}
  There is a unique left $\CH$-module structure on $K^{\Gbar}( \FNt
  \times \CB)$ with the property that for every $k$ in $K^{\Gbar}(
  \FNt \times \CB)$, $\lambda$ in $X(T)$, and simple reflection $s$ in
  $W$ we have
  \begin{itemize}
  \item[(a)] $-(T_s+1) \cdot k= l_s *k$ and
  \item[(b)] $E^\lambda \cdot k= e^\lambda *k$.
  \end{itemize}
\end{lemma}

Now the $\CH$-module and $K^{\Gbar}(Z)$-module structures on
$K^{\Gbar}( \FNt \times \CB)$ determine $A$-linear ring homomorphisms
$\phi_1\colon \CH\to \End_A \left( K^{\Gbar}( \FNt \times \CB)
\right)$ and $\phi_2\colon K^{\Gbar}(Z) \to \End_A\left( K^{\Gbar}(
  \FNt \times \CB) \right)$ respectively. It follows from Lemma
\ref{contain} that the image of $\phi_1$ is contained in the image of
$\phi_2$ and it follows from Proposition \ref{faithful} that $\phi_2$
is an injection.  Therefore, $\phi_2\inverse \circ \phi_1$ determines
an $A$-algebra homomorphism from $\CH$ to $K^{\Gbar}(Z)$ that we
denote by $\phi$.

The following theorem is proved in \cite[\S8]{lusztig:bases} using a
construction that goes back to \cite{kazhdanlusztig:langlands}. 

\begin{theorem}\label{iso}
  The $A$-algebra homomorphism $\phi\colon \CH\to K^{\Gbar}(Z)$ is an
  isomorphism and 
  \[
  \xymatrix{ Z(\CH) \ar@{^{(}->}[r] \ar[d]_{\cong} & A[X(T)]
    \ar@{^{(}->}[r] \ar[d]_{\cong}& \CH \ar[d]_{\cong} ^{\phi} \\
    K^{\Gbar} (\{ \operatorname{pt} \}) \ar@{^{(}->}[r] &
    K^{\Gbar}(Z_1) \ar@{^{(}->}[r] &K^{\Gbar}(Z)}
  \]
  is a commutative diagram of $A$-algebras and $A$-algebra
  homomorphisms.
\end{theorem}

In \cite[\S7.6]{chrissginzburg:representation} Chriss and Ginzburg
construct an isomorphism $\CH\cong K^{\Gbar}(Z)$ that satisfies the
conclusions of Theorem \ref{iso} using a variant of the ideas
above. 

Set $e= \sum_{w\in W} T_w$ in $\CH$. It is easy to check that there is
an $A$-module isomorphism $K^{\Gbar}(\FNt) \cong \CH e$ and hence an 
$A$-algebra isomorphism $\End_A( K^{\Gbar}( \FNt)) \cong \End_A( \CH
e)$. The convolution product construction can be used to define the
structure of a left $K^{\Gbar}(Z)$-module on $K^{\Gbar}(\FNt)$
\cite[\S5.4]{chrissginzburg:representation} and hence an$A$-algebra
homomorphism $K^{\Gbar}(Z) \to \End_A( K^{\Gbar}( \FNt))$. Similarly,
the left $\CH$-module structure on $\CH e$ defines an $A$-algebra
homomorphism $\CH \to \End_A( \CH e)$. Chriss and Ginzburg show that
the diagram
\[
\xymatrix{ \CH \ar[r]  & \End_A( \CH e) \ar[d]^{\cong} \\
K^{\Gbar}(Z) \ar[r] & \End_A( K^{\Gbar}( \FNt))}
\]
can be completed to a commutative square of $A$-algebras and that the
resulting $A$-algebra homomorphism $\CH\to K^{\Gbar}(Z)$ is an
isomorphism. We will see in \S\ref{s4.5} how this construction leads
to a conjectural description of the equivariant $K$-theory of the
generalized Steinberg varieties $\XPQ$.

\subsection{Irreducible representations of $\CH$, two-sided cells, and
  nilpotent orbits} \label{4.4} 

The isomorphism in Theorem \ref{iso} has been used by Kazhdan and
Lusztig \cite[\S 7]{kazhdanlusztig:langlands} to give a geometric
construction and parametrization of irreducible $\CH$-modules. Using
this construction, Lusztig \cite[\S 4]{lusztig:cellsIV} has found a
bijection between the set of two-sided Kazhdan-Lusztig cells in $W_e$
and the set of $G$-orbits in $\FN$. In order to describe this
bijection, as well as a conjectural description of two-sided ideals in
$K^{\Gbar}(Z)$ analogous to the decomposition of $H_{4n}(Z)$ given in
Proposition \ref{ideal}, we need to review the Kazhdan-Lusztig theory
of two-sided cells and Lusztig's based ring $J$.

The rules $v\mapsto v\inverse$ and $T_x\mapsto T_{x\inverse}
\inverse$, for $x$ in $W_e$, define a ring involution of $\CH$ denoted
by $h\mapsto \overline{h}$. The argument given by Kazhdan and Lusztig
in the proof of \cite[Theorem 1.1]{kazhdanlusztig:coxeter} applies to
$\CH$ and shows that there is a unique basis, $\{\, c_y'\mid y\in
W_e\,\}$, of $\CH$ with the following properties: 
\begin{itemize}
\item[(1)] $\overline{c_y'} = c_y'$ for all $y$ in $W_e$; and 
\item[(2)] if we write $c_y'= v^{-\ell(y)}\sum _{x\in W_e} P_{x,y}
  c_x'$, then $P_{y,y}=1$, $P_{x,y}=0$ unless $x\leq y$, and $P_{x,y}$
  is a polynomial in $v^2$ with degree (in $v$) at most $\ell(y)-
  \ell(x) -1$ when $x<y$.
\end{itemize}
The polynomials $P_{x,y}$ are called \emph{Kazhdan-Lusztig polynomials}.

For $x$ and $y$ in $W_e$, define $x\leq_{LR} y$ if there exists $h_1$
and $h_2$ in $\CH$ so that when $h_1 c_y' h_2$ is expressed as a
linear combination of $c_z'$, the coefficient of $c_x'$ is non-zero.
It follows from the results in \cite[\S 1]{kazhdanlusztig:coxeter}
that $\leq_{LR}$ is a preorder on $W_e$. The equivalence classes
determined by this preorder are \emph{two-sided Kazhdan-Lusztig
  cells.}

Suppose that $\Omega$ is a two-sided cell in $W_e$ and $y$ is in
$W_e$. Define $y\leq_{LR} \Omega$ if there is a $y'$ in $\Omega$ with
$y\leq_{LR} y'$. Then by construction, the span of $\{\, c_y'\mid
y\leq_{LR} \Omega \,\}$ is a two-sided ideal in $\CH$. We denote this
two-sided ideal by $\CH_{\overline{\Omega}}$.

The two sided ideals $\CH_{\overline{\Omega}}$ define a filtration of
$\CH$. In \cite[\S 2]{lusztig:cellsII}, Lusztig has defined a ring $J$
which after extending scalars is isomorphic to $\CH$, but for which
the two-sided cells index a decomposition into orthogonal two-sided
ideals, rather than a filtration by two-sided ideals.

For $x$, $y$, and $z$ in $W_e$, define $h_{x,y,z}$ in $A$ by $c_x'
c_y'= \sum_{z \in W_e} h_{x,y,z} c_z'$. Next, define $a(z)$ to be the
least non-negative integer $i$ with the property that $v^i h_{x,y,z}$
is in $\BBZ[v]$ for all $x$ and $y$. It is shown in \cite[\S
7]{lusztig:cellsI} that $a(z)\leq n$. Finally, define $\gamma_{x,y,z}$
to be the constant term of $v^{a(z)} h_{x,y,z}$.

Now let $J$ be the free abelian group with basis $\{\, j_y\mid y\in
W_e\,\}$ and define a binary operation on $J$ by $j_x* j_y =\sum_{z\in
  W_e} \gamma_{x,y,z} j_{z}$. For a two-sided cell $\Omega$ in $W_e$,
define $J_{\Omega}$ to be the span of $\{\, j_y\mid y\in \Omega \,\}$.
In \cite[\S 2]{lusztig:cellsII}, Lusztig proved that there are only
finitely many two-sided cells in $W_e$ and derived the following
properties of $(J, *)$:
\begin{itemize}
\item[(1)] $(J,*)$ is an associative ring with identity.
\item[(2)] $J_\Omega$ is a two-sided ideal in $J$ and $(J_\Omega, *)$
  is a ring with identity.
\item[(3)] $J \cong \oplus_{\Omega} J_\Omega$ (sum over all two-sided
  cells $\Omega$ in $W_e$).
\item[(4)] There is a homomorphism of $A$-algebras, $\CH\to J\otimes
  A$.
\end{itemize}

Returning to geometry, recall that $\CU$ denotes the set of unipotent
elements in $G$ and that $\CB_u= \{\, B'\in \CB\mid u\in B'\,\}$ for
$u$ in $\CU$.

Suppose $u$ is in $\CU$, $s$ in $G$ is semisimple, and $u$ and $s$
commute. Let $\langle s\rangle$ denote the smallest closed,
diagonalizable subgroup of $G$ containing $s$ and set $\sbar =\langle
s\rangle \times \BBC^*$. In \cite[\S2]{lusztig:cellsIV}, Lusztig
defines an action of $\sbar$ on $\CB_u$ using a homomorphism
$\SL_2(\BBC)\to G$ corresponding to $u$.  Define
\[
A_{\BBC}= A\otimes \BBC, \quad \CH_{\BBC}= \CH \otimes_A A_{\BBC},
\quad \text{and} \quad \CK_{u,s}= \left( K^{\sbar}( \CB_u) \otimes
  \BBC \right) \otimes_{R( \sbar)\otimes \BBC} A_{\BBC}.
\]
In \cite[\S2]{lusztig:cellsIV}, Lusztig defines commuting actions of
$\CH_{\BBC}$ and $C(us)$ on $\CK_{u,s}$. For an irreducible
representation $\rho$ of $C(us)$, let $\CK_{u,s, \rho}$ denote the
$\rho$-isotopic component of $\CK_{u,s}$, so $\CK_{u,s, \rho}$ is an
$\CH_{\BBC}$-module.  The next result is proved in \cite[Theorem
4.2]{lusztig:cellsIV}.

\begin{theorem}
  Suppose $u$ and $s$ are as above and that $\rho$ is an irreducible
  representation of $C(us)$ such that $\CK_{u,s, \rho}\ne 0$. Then,
  up to isomorphism, there is a unique simple $J$-module, $E$, with
  the property that when $E\otimes_{\BBC[v,v\inverse]} \BBC(v)$ is
  considered as an $\CH_{\BBC} \otimes_{\BBC[v,v\inverse]}
  \BBC(v)$-module, via the homomorphism $\CH\to J\otimes A$, then
  $E\otimes_{ \BBC[v,v\inverse]} \BBC(v) \cong \CK_{u,s,\rho}
  \otimes_{\BBC[v,v\inverse]} \BBC(v)$.
\end{theorem}

Given $u$, $s$, and $\rho$ as in the theorem, let $E(u,s,\rho)$ denote
the corresponding simple $J$-module. Since $J\cong \oplus_{\Omega}
J_\Omega$ and $E(u,s, \rho)$ is simple, there is a unique two-sided
cell $\Omega(u,s,\rho)$ with the property that $J_{\Omega(u,s,\rho)}
E(u, s, \rho) \ne0$. The main result in \cite[Theorem
4.8]{lusztig:cellsIV} is the next theorem.

\begin{theorem}
  With the notation as above, the two-sided cell $\Omega(u,s,\rho)$
  depends only on the $G$-conjugacy class of $u$. Moreover, the rule
  $(u,s,\rho) \mapsto \Omega(u,s,\rho)$ determines a well-defined
  bijection between the set of unipotent conjugacy classes in $G$ and
  the set of two-sided cells in $W_e$. This bijection has the property
  that $a(z)= \dim \CB_u$ for any $z$ in $\Omega(u,s,\rho)$.
\end{theorem}

Using a Springer isomorphism $\CU\cong \FN$ we obtain the following
corollary.

\begin{corollary}\label{bij}
  There is a bijection between the set of nilpotent $G$-orbits in
  $\FN$ and the set of two-sided cells of $W_e$ with the property that
  if $x$ is in $\FN$ and $\Omega$ is the two-sided cell corresponding
  to the $G$-orbit $G\cdot x$, then $a(z)= \dim \CB_x$ for every $z$
  in $\Omega$.
\end{corollary}

We can now work out some examples. Let $\Omega_1$ denote the two-sided
Kazhdan-Lusztig cell corresponding to the regular nilpotent orbit.
Then $a(z)=0$ for $z$ in $\Omega_1$ and $\Omega_1$ is the unique
two-sided Kazhdan-Lusztig cell on which the $a$-function takes the
value $0$. Let $1$ denote the identity element in $W_e$. Then it
follows immediately from the definitions that $\{1\}$ is a two-sided
cell and that $a(1)=0$. Therefore, $\Omega_1=\{1\}$.

At the other extreme, let $\Omega_0$ denote the two-sided
Kazhdan-Lusztig cell corresponding to the nilpotent orbit $\{0\}$.
Then $a(z)=n$ for $z$ in $\Omega_0$ and $\Omega_0$ is the unique
two-sided Kazhdan-Lusztig cell on which the $a$-function takes the
value $n$.  Shi \cite{shi:two} has shown that
\[
\Omega_0 = \{\, y\in W_e\mid a(y)=n\,\}= \{\, y_1w_0y_2\in W_e\mid
\ell(y_1 w_0y_2)= \ell(y_1) +\ell(w_0) +\ell(y_2) \,\}.
\]

The relation $\leq_{LR}$ determines a partial order on the set of
two-sided Kazhdan-Lusztig cells and one of the important properties of
Lusztig's $a$ function is that $a(y_1)\leq a(y_2)$ whenever
$y_2\leq_{LR} y_1$ (see \cite[Theorem
5.4]{lusztig:cellsI}). Therefore, $\Omega_1$ is the unique maximal
two-sided cell and $\Omega_0$ is the unique minimal two-sided cell. It
follows that $\CH_{\overline{ \Omega_1}}= \CH$ and that
$\CH_{\overline{ \Omega_0}}$ is the span of $\{\, c_y'\mid y\in
\Omega_0\,\}$.

Summarizing, we have seen that $\CH$ is filtered by the two sided
ideals $\CH_{\overline{\Omega}}$, where $\Omega$ runs over the set of
two-sided Kazhdan-Lusztig cells in $W_e$, and that there is a
bijection between the set of two-sided cells in $W_e$ and the set of
nilpotent orbits n $\FN$.

Now suppose that $\FO$ is a nilpotent orbit and recall the subvariety
$Z_{\overline{\FO}}$ of $Z$ defined in \S\ref{s3.5}. Let
$i_{\overline{\FO}}\colon Z_{\overline{\FO}} \to Z$ denote the
inclusion. There are direct image maps, $(i_{ \overline{ \FO}})_*$ in
Borel-Moore homology and in $K$-theory. It follows from the
convolution construction that the images of these maps are two-sided
ideals in $H_*(Z)$ and $K^{\Gbar}(Z)$ respectively. In \S\ref{s3.5} we
described the image of $(i_{\overline{\FO}})_* \colon H_{4n}( Z_{
  \overline{\FO}}) \to H_{4n}(Z)$, a two-sided ideal in $H_{4n}(Z)$.

The argument in \cite[\S5]{kazhdanlusztig:langlands} shows that
$(i_{\overline{\FO}})_* \otimes \id \colon K^{\Gbar}( Z_{ \overline{
    \FO}}) \otimes \BBQ \to K^{\Gbar}(Z) \otimes \BBQ$ is injective.
In contrast, $(i_{\overline{\FO}})_* \colon H_{j}( Z_{ \overline{
    \FO}}) \to H_{j}(Z)$ is an injection when $j=4n$, but fails to be
an injection in general. For example, taking $\overline{ \FO}=\FO=
\{0\}$, we have that $Z_{\{0\}}= \{0\}\times \CB\times \CB$ and $\dim
H_*(Z_{\{0\} })= \dim H_*(Z) =|W^2|$. However, $\dim H_{4n}(
Z_{\{0\}})=1$ and $H_{4n}(Z)=|W|$ and so $(i_{\{0\}})_* \colon
H_j(Z_{\{0\}}) \to H_j(Z)$ cannot be an injection for all $j$.

Define $\CI_{\overline{\FO}}$ to be the image of $(i_{\overline{
    \FO}})_* \colon K^{\Gbar}( Z_{ \overline{ \FO}}) \to
K^{\Gbar}(Z)$, a two-sided ideal in $K^{\Gbar}(Z)$. There is an
intriguing conjectural description of the image of $\CI_{ \overline{
    \FO}}$ under the isomorphism $K^{\Gbar}(Z) \cong \CH$ due to
Ginzburg \cite{ginzburg:geometrical} that ties together all the themes
in this subsection.

\begin{conjecture}
  Suppose that $\FO$ is a $G$-orbit in $\FN$ and $\Omega$ is the
  two-sided cell in $W_e$ corresponding to $\FO$ as in Corollary
  \ref{bij}. Then $\phi( \CI_{ \overline{ \FO}}) = \CH_{\overline{
      \Omega}}$, where $\phi\colon K^{\Gbar}(Z)\xrightarrow{\cong}
  \CH$ is the isomorphism in Theorem \ref{iso}.
\end{conjecture}

This conjecture has been proved when $G$ has type $A_l$ by Tanisaki
and Xi \cite{tanisakixi:kazhdan}. Xi has recently shown that the
conjecture is true after extending scalars to $\BBQ$
(\cite{xi:kazhdan}).

As a first example, consider the case of the regular nilpotent orbit
and the corresponding two-sided cell $\Omega_1$. Then $\overline{\FO}=
\FN$, $\CI_{\FN}= K^{\Gbar}(Z)$ and $\CH_{\overline{ \Omega_1}}= \CH$.
Thus the conjecture is easily seen to be true in this case. 

For a more interesting example, consider the case of the zero
nilpotent orbit. Then $Z_{\{0\}} = \{0\} \times \CB\times \CB$. The
corresponding two-sided cell, $\Omega_0$, has been described above and
we have seen that $\CH_{\overline{ \Omega_0}}$ is the span of $\{\,
c_y'\mid y\in \Omega_0\,\}$.

It is easy to check that $P_{w, w_0}=1$ for every $w$ in $W$ and thus
$c_{w_0}'= v^{-n} \sum_{w\in W} T_w= v^{-n}e$, where $e$ is as in
\S\ref{s4.3}. Let $\CH c_{w_0}' \CH$ denote the two sided ideal
generated by $c_{w_0}'$. In \cite{xi:representations}, Xi has proved
the following theorem.

\begin{theorem}\label{i*}
  With the notation as above we have
  \[
  \phi\left( \CI_{ \overline{ \{0\}}} \right)= \CH c_{w_0}' \CH =
  \CH_{\overline{\Omega_0}}.
  \]
\end{theorem}

\subsection{Equivariant $K$-theory of generalized Steinberg
  varieties}\label{s4.5}

Suppose $\CP$ and $\CQ$ are conjugacy classes of parabolic subgroups
of $G$ and recall the generalized Steinberg varieties $\XPQ$ and
$\XPQoo$, and the maps $\eta\colon Z\to \XPQ$ and $\eta_1\colon \ZPQ=
\eta\inverse( \XPQoo) \to \XPQoo$ from \S\ref{s2.4}. We have a
cartesian square of proper morphisms
\begin{equation}
  \label{diag}
  \xymatrix{
    \ZPQ \ar[r]_{k} \ar[d]_{\eta_1} & Z\ar[d]^{\eta} \\
    \XPQoo \ar[r]_{k_1} & \XPQ}
\end{equation}
where $k$ and $k_1$ are the inclusions. 

The morphism $\eta_1$ is smooth and so there is a pullback map in
equivariant $K$-theory, $\eta_1^*\colon K^{\Gbar}(\XPQoo) \to
K^{\Gbar}(\ZPQ)$. We can describe the $R(\Gbar)$-module structure of
$K^{\Gbar}(\ZPQ)$ and $K^{\Gbar}(\XPQoo)$ using the argument in
\cite[7.15]{lusztig:bases} together with a stronger concentration
theorem due to Thomason \cite[\S2]{thomason:formule}.

\begin{theorem}\label{kthy}
  The homomorphisms $\eta_1^*\colon K^{\Gbar}(\XPQoo) \to
  K^{\Gbar}(\ZPQ)$ and $k_*\colon K^{\Gbar} (\ZPQ) \to K^{\Gbar}(Z)$
  are injective. Moreover, $K^{\Gbar} (\XPQoo)$ is a free
  $R(\Gbar)$-module with rank $|W|^2 / |W_P| |W_Q|$ and $K^{\Gbar}
  (\ZPQ)$ is a free $R(\Gbar)$-module with rank $|W|^2$.
\end{theorem}

The Cellular Fibration Lemma of Chriss and Ginzburg
\cite[6.2.7]{chrissginzburg:representation} can be used to describe
the $R(\Gbar)$-module structure of $K^{\Gbar}(\XPQ)$ when $\CP=\CB$ or
$\CQ=\CB$.

\begin{proposition}
  The equivariant $K$-group $K^{\Gbar}( X^{\CB, \CQ})$ is a free
  $R(\Gbar)$-module with rank $|W|^2 / |W_Q|$.
\end{proposition}

We expect that $K^{\Gbar}(\XPQ)$ is a free $R(\Gbar)$-module with rank
$|W|^2 / |W_P| |W_Q$ for arbitrary $\CP$ and $\CQ$. We make a more
general conjecture about $K^{\Gbar}(\XPQ)$ after first considering an
example in which everything has been explicitly computed. 

Consider the very special case when $\CP=\CQ=\{G\}$. In this case the
spaces in (\ref{diag}) are well-known:
\[
X^{\{G\}, \{G\}} _{0,0}\equiv \{0\},\quad Z^{\{G\}, \{G\}} =
\overline{ Z_{w_0}}= Z_{\{0\}} \cong \CB \times \CB,\quad \text{and}
\quad X^{\{G\}, \{G\}}\equiv \FN.
\]
Also, $\eta\colon Z\to X^{\{G\}, \{G\}}$ may be identified with
$\mu_z\colon Z\to \FN$ and $k\colon Z^{\{G\}, \{G\}}\to Z$ may be
identified with the closed embedding $\CB\times \CB \to Z$ by $(B',
B'')\mapsto (0, B', B'')$ and so (\ref{diag}) becomes
\[
\xymatrix{ Z^{\{G\},\{G\}} =Z_{\{0\}} \ar[r]_-{i_{\{0\}}} \ar[d] &
  Z\ar[d]^{\mu_z} \\
  X^{\{G\},\{G\}}_{0,0}= \{0\} \ar[r] & \FN \cong X^{\{G\},\{G\}}.}
\]
The image of $(i_{\{0\}})_* \colon K^{\Gbar} (Z_{\{0\}}) \to K^{\Gbar}
(Z)$ is $\CI_{ \overline{ \{0\}}}$ and we saw in Theorem \ref{i*} that
$\CI_{ \overline{ \{0\}}} \cong \CH c_{w_0}' \CH =
\CH_{\overline{\Omega_0}}$.

Ostrik \cite{ostrik:equivariant} has described the map $(\mu_z)_*
\colon K^{\Gbar}(Z)\to K^{\Gbar}(X^{\{G\},\{G\}})$.  Recall that $W_e=
X(T) \rtimes W$. Because the fundamental Weyl chamber is a fundamental
domain for the action of $W$ on $X(T)\otimes \BBR$, it follows that
each $(W,W)$-double coset in $W_e$ contains a unique element in $X^+$.
Also, each $(W,W)$-double coset in $W_e$ contains a unique element
with minimal length. For $\lambda$ in $X^+$ we let $m_\lambda$ denote
the element with minimal length in the double coset $W\lambda W$.

\begin{theorem}
  For $x$ in $W_e$, $(\mu_z)_*(c_x')=0$ unless $x=m_\lambda$ for some
  $\lambda$ in $X^+$. Moreover, the map $(\mu_z)_*\colon K^{\Gbar}(Z)
  \to K^{\Gbar}( X^{\{G\},\{G\}})$ is surjective and $\{\,
  (\mu_z)_*(c_{m_\lambda}') \mid \lambda\in X^+\,\}$ is an $A$-basis
  of $K^{\Gbar}( X^{\{G\}, \{G\}})$.
\end{theorem}

Notice that the theorem is the $K$-theoretic analog of Theorem
\ref{ave} in the very special case we are considering. 

To prove this result, Ostrik uses the description of $Z$ as a fibred
product and the two corresponding factorizations of $\mu_z$:
\begin{equation}
  \label{cdK}
  \xymatrix{ Z =\FNt\times_{\FN} \FNt \ar[r] \ar[d] &
    X^{\CB, \{G\}} \cong \FNt  \ar[d] \\
    \FNt \cong X^{\{G\},\CB} \ar[r] & X^{\{G\},\{G\}} \cong \FN.}  
\end{equation}
It follows from the construction of the isomorphism $K^{\Gbar}(Z)
\cong \CH$ given by Chriss and Ginzburg
\cite[\S7.6]{chrissginzburg:representation} (see the end of
\S\ref{s4.3}) that after applying the functor $K^{\Gbar}$ to
(\ref{cdK}) the resulting commutative diagram of equivariant
$K$-groups may be identified with the following commutative diagram
subspaces of $\CH$:
\begin{equation} \label{cdH}
  \xymatrix{ \CH \ar[r] \ar[d] &
    \CH c_{w_0}'  \ar[d] \\
    c_{w_0}' \CH \ar[r] & c_{w_0}' \CH c_{w_0}' }
\end{equation}
where the maps are given by the appropriate right or left
multiplication by $c_{w_0}'$.

We conclude with a conjecture describing $K^{\Gbar}(\XPQ)$ for
arbitrary $\CP$ and $\CQ$. Recall from \S\ref{s3.7} that $\XPQ\cong
\FNt_{\CP} \times_{\FN} \FNt^{\CQ}$. The projection $\mu\colon \FNt\to
\FN$ factors as $\FNt \xrightarrow{\eta^{\CP}} \FNt^{\CP}
\xrightarrow{\xi^{\CP}} \FN$ where $\eta^{\CP}(x, gBg\inverse)= (x,
gPg\inverse)$ and $\xi^{\CP}(x, gPg\inverse)= x$. Using this
factorization, we may expand diagram (\ref{cdK}) to a $3\times 3$
diagram with $\XPQ$ in the center:
\begin{equation}
  \label{3x3}
  \xymatrix{ Z\ar[r] \ar[d] &X^{\CB, \CQ} \ar[r] \ar[d] & \FNt
    \ar[d]\\
    X^{\CP, \CB}\ar[r] \ar[d] &X^{\CP, \CQ} \ar[r] \ar[d] & \FNt^{\CP}
    \ar[d]\\
    \FNt \ar[r] & \FNt^{\CQ} \ar[r] & \FN. }  
\end{equation}

Let $w_P$ and $w_Q$ denote the longest elements in $W_P$ and $W_Q$
respectively. Comparing (\ref{cdK}), (\ref{cdH}), and (\ref{3x3}), we
make the following conjecture. This conjecture is a $K$-theoretic
analog of (\ref{X}) and Conjecture \ref{conj}.

\begin{conjecture}\label{conj2}
  With the notation above, $K^{\Gbar}(\XPQ) \cong c_{w_P}' \CH
  c_{w_Q}'$.
\end{conjecture}

If the conjecture is true, then after applying the functor $K^{\Gbar}$
to (\ref{3x3}) the resulting commutative diagram of equivariant
$K$-groups may be identified with the following commutative diagram of
subspaces of $\CH$:
\begin{equation}
  \xymatrix{ \CH \ar[r] \ar[d] & \CH c_{w_Q}' \ar[r] \ar[d] &\CH
    c_{w_0}'  \ar[d] \\ 
    c_{w_P}' \CH \ar[r] \ar[d] & c_{w_P}' \CH c_{w_Q}' \ar[r] \ar[d]&
    c_{w_P}' \CH  c_{w_0}'  \ar[d] \\ 
    c_{w_0}' \CH \ar[r] & c_{w_0}' \CH c_{w_Q}' \ar[r] & c_{w_0}' \CH
    c_{w_0}'. } 
\end{equation}



\bigskip

\bibliographystyle{amsalpha}

\begin{thebibliography}{BBM89}

\bibitem[BB81]{beilinsonbernstein:localisation} A.~Beilinson and
  J.~Bernstein, \emph{Localisation de {$g$}-modules}, C. R.
  Acad. Sci. Paris S\'er. I Math. \textbf{292} (1981), no.~1, 15--18.

\bibitem[BB82]{borhobrylinski:differentialI} W.~Borho and
  J.L. Brylinski, \emph{Differential operators on homogeneous
    spaces. {I} {I}rreducibility of the associated variety for
    annihilators of induced modules.}, Invent. Math. \textbf{69}
  (1982), 437--476.

\bibitem[BB85]{borhobrylinski:differentialIII} \bysame,
  \emph{Differential operators on homogeneous spaces {III};
    {C}haracteristic varieties of {H}arish-{C}handra modules and of
    primitive ideals}, Invent. Math. \textbf{80} (1985), 1--68.

\bibitem[BBD82]{beilinsonbernsteindeligne:faisceaux} A.~Beilinson,
  J.~Bernstein, and P.~Deligne, \emph{Faisceaux pervers}, Analysis and
  topology on singular spaces, I (Luminy, 1981), Ast\'erisque,
  vol. 100, Soc. Math. France, Paris, 1982, pp.~5--171.

\bibitem[BBM89]{borhobrylinskimacpherson:nilpotent} W.~Borho,
  J.L. Brylinski, and R.~MacPherson, \emph{Nilpotent {O}rbits,
    {P}rimitive {I}deals, and {C}haracteristic {C}lasses}, Progress in
  Mathematics, vol.~78, Birkh{\"a}user, Boston, MA, 1989.

\bibitem[BG80]{bernsteingelfand:tensor} J.~Bernstein and S.I. Gelfand,
  \emph{Tensor products of finite- and infinite-dimensional
    representations of semisimple {L}ie algebras}, Compositio
  Math. \textbf{41} (1980), no.~2, 245--285.

\bibitem[BM81]{borhomacpherson:weyl} W.~Borho and R.~MacPherson,
  \emph{Repr{\'e}sentations des groupes de {W}eyl et homologie
    d'intersection pour les vari{\'e}t{\'e}s nilpotentes}, C. R. Acad.
  Sci. Paris \textbf{292} (1981), no.~15, 707--710.

\bibitem[BM83]{borhomacpherson:partial} \bysame, \emph{Partial
    resolutions of nilpotent varieties. {A}nalysis and topology on
    singular spaces, {II}, {III} ({L}uminy, 1981)}, Ast\'erisque
  \textbf{101} (1983), 23--74.

\bibitem[Bor84]{borel:intersection} A.~Borel (ed.), \emph{Intersection
    cohomology ({B}ern, 1983)}, Progr. Math., vol.~50, Birkh{\"a}user,
  Boston, 1984.

\bibitem[Car85]{carter:finite} R.~Carter, \emph{Finite groups of {L}ie
    type}, Pure and Applied Mathematics (New York), John Wiley \& Sons
  Inc., New York, 1985, Conjugacy classes and complex characters, A
  Wiley-Interscience Publication.

\bibitem[CG97]{chrissginzburg:representation} N.~Chriss and
  V.~Ginzburg, \emph{Representation theory and complex geometry},
  Birkh{\"{a}}user, Boston, 1997.

\bibitem[Dim04]{dimca:sheaves} A.~Dimca, \emph{Sheaves in topology},
  Universitext, Springer-Verlag, Berlin, 2004.

\bibitem[Dou96]{douglass:involution} J.M. Douglass, \emph{An
    involution of the variety of flags fixed by a unipotent linear
    transformation}, Adv. Appl. Math. \textbf{17} (1996), 357--379.

\bibitem[DR04]{douglassroehrle:geometry} J.M. Douglass and
  G.~R{\"o}hrle, \emph{The geometry of generalized {S}teinberg
    varieties}, Adv. Math. \textbf{187} (2004), no.~2, 396--416.

\bibitem[DR08a]{douglassroehrle:homology} \bysame,
  \emph{{B}orel-{M}oore homology of generalized {S}teinberg
    varieties}, Trans. AMS (2008), to appear.

\bibitem[DR08b]{douglassroehrle:coinvariant} \bysame, \emph{{H}omology
    of the {S}teinberg variety and {W}eyl group coinvariants},
  arXiv:0704.1717v1.

\bibitem[FM81]{fultonmacpherson:categorical} W.~Fulton and
  R.~MacPherson, \emph{Categorical framework for the study of singular
    spaces}, Mem. Amer. Math. Soc. \textbf{31} (1981), no.~243.

\bibitem[Gin86]{ginzburg:G-modules} V.~Ginzburg, \emph{{${\mathfrak
        G}$}-modules, {S}pringer's representations and bivariant
    {C}hern classes}, Adv. in Math. \textbf{61} (1986), no.~1, 1--48.

\bibitem[Gin87]{ginzburg:geometrical} \bysame, \emph{Geometrical
    aspects of representation theory}, Proceedings of the
  International Congress of Mathematicians, Vol. 1, 2 (Berkeley,
  Calif., 1986) (Providence, RI), Amer. Math. Soc., 1987,
  pp.~840--848.

\bibitem[GM83]{goreskymacpherson:intersectionII} M.~Goresky and
  R.~MacPherson, \emph{Intersection homology. {II}}, Invent. Math.
  \textbf{72} (1983), no.~1, 77--129.

\bibitem[HJ05]{hinichjoseph:orbital} {V}. Hinich and {A}. Joseph,
  \emph{Orbital variety closures and the convolution product in
    {B}orel-{M}oore homology}, Selecta Math. (N.S.) \textbf{11}
  (2005), no.~1, 9--36.

\bibitem[Hot82]{hotta:springer} R.~Hotta, \emph{On {S}pringer's
    representations}, J. Fac. Sci. Univ. Tokyo Sect. IA Math. (1982),
  no.~3, 863--876.

\bibitem[Hot85]{hotta:local} \bysame, \emph{A local formula for
    {S}pringer's representation}, Algebraic groups and related topics
  ({K}yoto/{N}agoya, 1983), Adv. {S}tud. {P}ure {M}ath.,
  North-Holland, Amsterdam, 1985, pp.~127--138.

\bibitem[IM65]{iwahorimatsumoto:bruhat} N.~Iwahori and H.~Matsumoto,
  \emph{On some {B}ruhat decompositions and the structure of the
    {H}ecke ring of $p$-adic {C}hevalley groups}, Publ. Math IHES
  \textbf{25} (1965), 5--48.

\bibitem[Jos79]{joseph:dixmier} A.~Joseph, \emph{Dixmier's problem for
    {V}erma and principal series submodules}, J. London Math. Soc. (2)
  \textbf{20} (1979), no.~2, 193--204.

\bibitem[Jos84]{joseph:variety} \bysame, \emph{On the variety of a
    highest weight module}, J. Algebra \textbf{88} (1984), 238--278.

\bibitem[KL79]{kazhdanlusztig:coxeter} D.~Kazhdan and G.~Lusztig,
  \emph{Representations of {C}oxeter groups and {H}ecke algebras},
  Invent. Math. \textbf{53} (1979), 165--184.

\bibitem[KL80]{kazhdanlusztig:topological} \bysame, \emph{A
    topological approach to {S}pringer's representations}, Adv. in
  Math. \textbf{38} (1980), 222--228.

\bibitem[KL87]{kazhdanlusztig:langlands} \bysame, \emph{Proof of the
    {D}eligne-{L}anglands conjecture for {H}ecke algebras},
  Invent. Math. \textbf{87} (1987), 153--215.

\bibitem[KS90]{kashiwaraschapira:sheaves} M.~Kashiwara and
  P.~Schapira, \emph{Sheaves on manifolds}, Grundlehren der
  Mathematischen Wissenschaften [Fundamental Principles of
  Mathematical Sciences], vol. 292, Springer-Verlag, Berlin, 1990,
  With a chapter in French by Christian Houzel.

\bibitem[KT84]{kashiwaratanisaki:characteristic} M.~Kashiwara and
  T.~Tanisaki, \emph{The characteristic cycles of holonomic systems on
    a flag manifold related to the {W}eyl group algebra}, Invent.
  Math. \textbf{77} (1984), no.~1, 185--198.

\bibitem[Lus79]{lusztig:class} G.~Lusztig, \emph{A class of
    irreducible representations of a {W}eyl group},
  Nederl. Akad. Wetensch. Indag. Math. \textbf{41} (1979), no.~3,
  323--335.

\bibitem[Lus84]{lusztig:characters} \bysame, \emph{Characters of
    reductive groups over a finite field}, Annals of Mathematics
  Studies, vol. 107, Princeton University Press, Princeton, NJ, 1984.

\bibitem[Lus85]{lusztig:cellsI} \bysame, \emph{Cells in affine {W}eyl
    groups}, Algebraic groups and related topics (Kyoto/Nagoya, 1983),
  Adv. Stud. Pure Math., vol.~6, North-Holland, Amsterdam, 1985,
  pp.~255--287.

\bibitem[Lus87]{lusztig:cellsII} \bysame, \emph{Cells in affine {W}eyl
    groups. {II}}, J. Algebra \textbf{109} (1987), no.~2, 536--548.

\bibitem[Lus89a]{lusztig:affine} \bysame, \emph{Affine {H}ecke
    algebras and their graded version}, J. Amer.
  Math. Soc. \textbf{2} (1989), no.~3, 599--635.

\bibitem[Lus89b]{lusztig:cellsIV} \bysame, \emph{Cells in affine
    {W}eyl groups. {IV}}, J. Fac. Sci. Univ. Tokyo Sect. IA
  Math. \textbf{36} (1989), no.~2, 297--328.

\bibitem[Lus98]{lusztig:bases} \bysame, \emph{Bases in equivariant
    {$K$}-theory}, Represent. Theory \textbf{2} (1998), 298--369
  (electronic).

\bibitem[Ost00]{ostrik:equivariant} V.~Ostrik, \emph{On the
    equivariant {$K$}-theory of the nilpotent cone}, Represent. Theory
  \textbf{4} (2000), 296--305 (electronic).

\bibitem[Shi87]{shi:two} J.Y. Shi, \emph{A two-sided cell in an affine
    {W}eyl group}, J. London Math.  Soc. (2) \textbf{36} (1987),
  no.~3, 407--420.

\bibitem[Sho88]{shoji:geometry} T.~Shoji, \emph{Geometry of orbits and
    {S}pringer correspondence}, Ast\'erisque (1988), no.~168, 9,
  61--140, Orbites unipotentes et repr\'esentations, I.

\bibitem[Slo80]{slodowy:simple} P.~Slodowy, \emph{Simple singularities
    and simple algebraic groups}, Lecture Notes in Mathematics,
  vol. 815, Springer-Verlag, Berlin/Heidelberg/New~York, 1980.

\bibitem[Spa76]{spaltenstein:fixed} N.~Spaltenstein, \emph{The fixed
    point set of a unipotent transformation on the flag manifold},
  Proc. Kon. Nederl. Akad. Wetensch. \textbf{79} (1976), 452--458.

\bibitem[Spa82]{spaltenstein:classes} \bysame, \emph{Classes unipotent
    et sous-groupes de {B}orel}, Lecture Notes in Mathematics,
  vol. 946, Springer-Verlag, Berlin/Heidelberg/New~York, 1982.

\bibitem[Spr76]{springer:trigonometric} T.A. Springer,
  \emph{Trigonometric sums, {G}reen functions of finite groups and
    representations of {W}eyl groups}, Invent. Math. \textbf{36}
  (1976), 173--207.

\bibitem[Spr78]{springer:construction} \bysame, \emph{A construction
    of representations of {W}eyl groups}, Invent.  Math. \textbf{44}
  (1978), 279--293.

\bibitem[Spr98]{springer:algebraic} T.A. Springer, \emph{Linear
    algebraic groups}, second ed., Progress in Mathematics, vol.~9,
  Birkh\"auser Boston Inc., Boston, MA, 1998.

\bibitem[Ste76]{steinberg:desingularization} R.~Steinberg, \emph{On
    the desingularization of the unipotent variety}, Invent.
  Math. \textbf{36} (1976), 209--224.

\bibitem[Ste88]{steinberg:robinson} \bysame, \emph{An occurrence of
    the {R}obinson-{S}chensted correspondence}, J.  Algebra
  \textbf{113} (1988), 523--528.

\bibitem[Tho92]{thomason:formule} R.W. Thomason, \emph{Une formule de
    {L}efschetz en {$K$}-th{\'e}orie {\'e}quivariante alg{\'e}brique},
  Duke Math. J. \textbf{68} (1992), no.~3, 447--462.

\bibitem[TX06]{tanisakixi:kazhdan} T.~Tanisaki and N.~Xi,
  \emph{Kazhdan-{L}usztig basis and a geometric filtration of an
    affine {H}ecke algebra}, Nagoya Math. J. \textbf{182} (2006),
  285--311.

\bibitem[Xi94]{xi:representations} N.~Xi, \emph{Representations of
    affine {H}ecke algebras}, Lecture Notes in Mathematics, vol. 1587,
  Springer-Verlag, Berlin, 1994.

\bibitem[Xi08]{xi:kazhdan} \bysame, \emph{{K}azhdan-{L}usztig basis
    and a geometric filtration of an affine {H}ecke algebra, {II}},
  arXiv:0801.0472.

\end{thebibliography}

\providecommand{\bysame}{\leavevmode\hbox to3em{\hrulefill}\thinspace}
\providecommand{\MR}{\relax\ifhmode\unskip\space\fi MR }
\providecommand{\MRhref}[2]{%
  \href{http://www.ams.org/mathscinet-getitem?mr=#1}{#2}
}
\providecommand{\href}[2]{#2}


\end{document}